\documentclass[final]{amsart}

\usepackage{microtype}
\usepackage{booktabs}
\usepackage{cite}
\usepackage{enumitem}
\usepackage{lipsum}
\usepackage{amsfonts}
\usepackage{graphicx}
\usepackage{epstopdf}
\usepackage{algorithmic}
\usepackage{algorithm}
\usepackage{amsopn}
\usepackage{amssymb}
\usepackage{amsmath}
\usepackage{amsthm}
\usepackage[colorlinks=true]{hyperref}
\usepackage{todonotes}
\usepackage{tikz-cd}
\usepackage{diagbox}
\usepackage{multirow}
\ifpdf
  \DeclareGraphicsExtensions{.eps,.pdf,.png,.jpg}
\else
  \DeclareGraphicsExtensions{.eps}
\fi

\newtheorem{theorem}{Theorem}
\newtheorem{lemma}[theorem]{Lemma}

\theoremstyle{definition}

\newtheorem*{example}{Example}

\theoremstyle{remark}
\newtheorem{remark}[theorem]{Remark}

\def\th{\theta}

\def\F{\mathcal{F}}

\def\B{\mathcal{B}}
\def\F{\mathcal{F}}
\def\G{\mathcal{G}}
\def\M{\mathcal{M}}
\def\I{\mathcal{I}}
\def\T{\mathcal{T}}
\def\K{\mathcal{K}}

\def\bB{{\boldsymbol B}}
\def\bV{{\mathbf V}}

\let\on=\operatorname

\def\R{\mathbb{R}}
\def\S{\mathbb{S}}

\let\on=\operatorname

\let\mf=\mathfrak

\def\Imm{\on{Imm}}
\def\Emb{\on{Emb}}
\def\Diff{\on{Diff}}
\def\vol{\on{vol}}
\def\Null{\on{Null}}
\def\id{{\on{id}}}
\def\att{\mathit{att}}

\newcommand{\scp}[2]{{\left\langle {#1}\, , \, {#2}\right\rangle}}

\graphicspath{{./Figures/}}

\begin{document}

\title{Metric Registration of Curves and Surfaces using Optimal Control}

\author{Martin Bauer}
\address{Department of Mathematics\\Florida State University}
\email{bauer@math.fsu.edu}

\author{Nicolas Charon}
\address{Department of Applied Mathematics and Statistics\\Johns Hopkins University}
\email{charon@cis.jhu.edu}

%    Information for third author
\author{Laurent Younes}
\address{Department of Applied Mathematics and Statistics\\Johns Hopkins University}
\email{laurent.younes@jhu.edu}

\renewcommand{\shortauthors}{M. Bauer, N. Charon and L. Younes}



\maketitle

\begin{abstract}
    This chapter presents an overview of recent developments in the analysis of shapes such as curves and surfaces through Riemannian metrics. We show that several constructions of metrics on spaces of submanifolds can be unified through the prism of Riemannian submersions, with shape space metrics being induced from metrics defined on the top spaces. Computing the resulting Riemannian distances involves solving geodesic matching problems with boundary conditions. To deal efficiently with such variational problems, one can rely on an auxiliary family of "chordal" distances to simplify the treatment of boundary conditions, which we use to come up with a relaxed inexact formulation of the matching problem. This also allows to turn shape matching into optimal control problems and give a common framework to address them in practice. We then specify our analysis to the cases of intrinsic shape metrics defined using invariant Sobolev metrics on parametrized immersions, outer shape metrics induced from metrics on diffeomorphism groups of the ambient space and finally a recent hybrid model that combines those two approaches.       
\end{abstract}

\section{Introduction}
Curves and surfaces are natural representations of shapes in two or three dimensions and are mathematically special cases of submanifolds of the Euclidean space $\R^d$ for $d=2$ or $3$. While such submanifolds are conveniently described through parametrizations, or atlases, which are relatively simple mathematical objects (forming special classes of functions), shapes, as geometric objects, form equivalence classes of such functions modulo reparametrization, leading to quotient spaces that can be significantly more complex to study, especially in infinite dimensions. We will refer to these equivalence classes of parametrized submanifolds modulo reparametrization, as ``unparametrized submanifolds'', mainly for clarity since these objects are what is commonly referred to as, simply, submanifolds.

Constructing tractable comparison tools between such unparametrized sets is at the core of shape analysis and its many applications. Among the various approaches that were proposed, our goal is here to focus the review on a class of metrics defined through the introduction of Riemannian (or sub-Riemannian) metrics on spaces of parametrized submanifolds with suitable invariance properties, allowing these metrics, through mechanisms that will be described below, to ``project'' onto distances modulo reparametrization. The computation of these distances will always require the estimation of a geodesic path between them, resulting in a class of algorithms that we will refer to as ``metric registration''.\\

In our setting, parametrized manifolds will be represented as functions $q:M \to \R^d$, where $M$, the parameter space, is a  $p$-dimensional oriented compact manifold. Typical examples are $M = \S^{p}$, the $p$-dimensional unit sphere, or a regular closed subset of $\R^p$ with a smooth boundary (such as a closed interval, or a closed disc). Situations in which $M$ has no boundary are usually more amenable to theoretical analysis, but the with-boundary case also has important practical applications (e.g., for open curves). 

In order for $q$ to parametrize a $p$-dimensional submanifold of $\R^d$, one needs to assume that it provides an embedding, i.e., a one-to-one $C^1$ mapping from $M$ into $\R^d$, such that $q^{-1}: q(M) \to M$ is continuous and $dq(m)$ is a one-to-one linear mapping for all $m\in M$ (where $dq$ denotes the differential of $q$). The set of embeddings will be denoted $\Emb(M, \R^d)$. However, for several metrics described in this paper, it will be useful to relax the embedding constraint and assume that $q$ is an immersion instead, where an immersion is a $C^1$ mapping from $M$ into $\R^d$ such that $dq(m)$ is a one-to-one linear mapping for all $m\in M$ (so that an embedding is an immersion which is, in addition, a homeomorphism onto its image). We let $\Imm(M, \R^d)$ denote the space of immersions from $M$ to $\R^d$. 
The main reason for considering immersions rather than embeddings is that several of the metrics that will be introduced  have geodesic paths that remain immersions at all times, but not necessarily embeddings (even when looking for shortest paths between two embeddings).   

If $X$ is a space of $C^1$ functions from $M$ to $\R^d$, we will use the notation $\Emb_X$, $\Imm_X$ for the set of embeddings or immersions $q\in X$. We will typically use $X = C^r =  C^r(M, \R^d)$ (functions with at least $r$ continuous derivatives) with $r\geq 1$, or $X = H^s = H^s(M, \R^d)$ (functions with at least $s$ generalized derivatives which are all square integrable) with $s>p/2+1$, or $X=C^{\infty}= C^{\infty}(M, \R^d)$ (smooth functions).
\\

Let $\Diff(M)$ denote the space of orientation-preserving diffeomorphisms of $M$, i.e., the space of invertible $C^1$ mappings $\tau$ from $M$ onto itself with a $C^1$ inverse such that, for all $m\in M$, $d\tau(m)$ maps positively oriented bases of $T_mM$ onto positively oriented bases of $T_{\tau(m)}M$. 
Given $q\in \Imm(M, \R^d)$, we let $S_q = q(M)$ denote the image of $M$ by $q$. If $\tau\in \Diff(M)$, the mapping $q\cdot \tau = q\circ \tau$ is also an immersion and $S_{q\cdot \tau} = S_q$. One says that $\tau$ is a change of parameter, or reparametrization, and the unparametrized surface associated with $q$ is the equivalent class
\[
[q] = q\cdot \Diff(M) = \{q\circ \tau, \tau\in\Diff(M)\}.
\]
% \todomartin{Maybe we should write 
% $q\cdot \Diff(M)$ to emphasize that it is a right action?}
% \todonicolas{We should decide on a convention for that part. I think we said we would write $\psi$ for elements of $\Diff(M)$. Then we either take $q \circ \psi^{-1}$ with a left action or $q\circ \psi$ with a right action. I believe the second option is maybe more common in Martin's papers.}
% \todolaurent{I did not want to use $\psi$ to keep the notation available for the outer deformations. Also, I think it will be easier for the reader if we only have left actions.}
Notice that, when considering immersions that belong to a given class $X$, the space $\Diff(M)$ must also be restricted to diffeomorphisms that preserve $X$, i.e., such that $q\circ \tau\in X$ for all $q\in X$.
\\

As indicated earlier, it will be relatively ``easy'' to equip $\Imm_X(M, \R^d)$ with a Riemannian metric and various approaches will be described later in this chapter. We will then use a standard mechanism --Riemannian submersions, described in the next section-- to define, through the mapping $q\mapsto [q]$, a distance between unparametrized shapes. More precisely, this chapter is organized as follows. Following this introduction, Section \ref{sec:submersion} describes the general concept of Riemannian submersion and how it can be utilized to define metrics in spaces acted upon by deformation groups. Section \ref{sec:control} summarizes basic results in optimal control, that will be instrumental in the computation of the geodesic distances. Section \ref{sec:chordal} reviews a special class of metrics in which shapes spaces are embedded in suitably chosen Hilbert spaces and compared using the associated Hilbert norm, yielding ``chordal metrics'' between shapes, to distinguish them from Riemannian metrics directly defined on shape spaces, which are our main focus. Chordal metrics will here mainly serve as auxiliary terms in order to relax variational problems associated to the computation of geodesics into computationally feasible optimal control problems. Section \ref{sec:intrinsic_metrics} introduces a first class of metrics, called intrinsic, because they derive from standard Sobolev metrics on the space of immersions. Section \ref{sec:LDDMM} describes ``outer-deformation metrics'' that are derived from the action of groups of diffeomorphisms of the ambient space on immersions. Finally, Section \ref{sec:hybrid} describes some attempts at combining the two approaches, leading to ``hybrid metrics.'' These last three sections will also include some experimental results and a discussion of numerical methods. We will conclude the chapter by a brief discussion in Section \ref{sec:conclusion}. 

\section{Building Metrics via Submersions}
\label{sec:submersion}

The following discussion describes some metric constructions involving Riemannian submersions. They can be rigorously defined in finite dimensions, but also serve as motivational foundations for similar constructions in infinite dimensions (which will be our main focus), even though not all of the properties listed below may generalize, or they may require additional assumptions to hold.

If $\F$ and and $\B$ are differential manifolds, with $\dim(\B) \leq \dim(\F)$, a mapping $\pi$ from $\F$ onto $\B$ is a submersion if $d\pi(f)$ is onto from $T_f\F$ to $T_{\pi(f)}\B$ for all $f\in \F$.
For such an immersion, one defines
\[
V_f = \Null(d\pi(f)) = \{\eta \in T_f\F: d\pi(f)\eta=0\},
\]
the vertical space at $f$ for the submersion $\pi$, which is also the tangent space at $f$ to the submanifold $\pi^{-1}(b)$, where $b=\pi(f)$.
%\todomartin{I would prefer vertical and horizontal bundle, instead of space}
If $\F$ is a Riemannian manifold, one can also define, for $f\in \F$, the horizontal space
\[
H_f = V_f^\perp = \{\eta\in T_f\F: g^\F_f(\eta,\zeta) =0\text{ for all } \zeta\in V_f\},
\]
where $g^\F$ denotes the Riemannian metric on $\F$, so that $g^\F_f(\eta, \zeta)$ is the Riemannian inner product between $\eta$ and $\zeta$ in $T_f\F$. Moreover, the restriction of $d\pi(f)$ to $H_f$ is an isomorphism between $H_f$ and $T_{b}\B$ where $b= \pi(f)$. If $\B$ is also Riemannian, then one says that $\pi$ is a {\em Riemannian submersion} if and only if $d\pi(f)$ is an isometry from $H_f$ to $T_b\B$, i.e., for all $\eta,\eta'\in H_f$:
\begin{equation}
\label{eq:r.sub}
g^\B_b(d\pi(f)\eta, d\pi(f)\eta') = g^\F_f(\eta,\eta').
\end{equation}

If $\pi(f)=\pi(f')=b$, this implies that the horizontal spaces $H_f$ and $H_{f'}$ are isometric, via the transformations $\psi_{f,f'} = d\pi(f')^{-1}\circ d\pi(f): H_f\to H_{f'}$. Conversely, if the metric on $\F$ is such that all transformations $\psi_{f,f'}$ are isometries whenever $\pi(f) = \pi(f')$, then one can equip $\B$ with a uniquely defined Riemannian metric, through \eqref{eq:r.sub}, such that $\pi$ is a Riemannian submersion. 
%\todomartin{This is in infinte dimensions not always guaranteed: the horizontal bundle does not necessarily exist. But i am ok if we want to ignore this technicality} 
This is the construction that we will use in the next sections, in which $\F$ will be a set on which Riemannian metrics can be easily defined, with properties ensuring that horizontal spaces are isometric, resulting in a metric on $\B$.\\

Let $\pi$ be a Riemannian submersion and $h\in T_b\B$ for some $b\in \B$. If $f\in \pi^{-1}(b)$, denote by $h^f$ the unique vector in $H_f$  such that $d\pi(f)h^f = h$, so that $g^\B_b(h,h) = g^\F_b(h^f,h^f) $. Any vector $\eta \in d\pi(f)^{-1}h$ can be uniquely decomposed as $\eta = h^f + \zeta$, where $\zeta\in V_f$, so that  the metric on $\B$ can also be derived from the one on $\F$ with
\begin{equation}
    \label{eq:r.sub.2}
g^\B_b(h,h) = \min\{g^\F_f(\eta, \eta): d\pi(f)\eta = h\}
\end{equation}
for any $f$ such that $\pi(f)=b$. This immediately implies that minimizing geodesics in $\B$, which are functions $t\mapsto b(t)$ such that $b(0)=b_0$ and $b(1)=b_1$ for given $b_0, b_1$ and that minimize
\[
\int_0^1 g_{b(t)}^\B(\dot b(t), \dot b(t)) dt
\]
can also be defined as $b(t) = \pi(f(t))$ where $f(\cdot)$ minimizes 
\begin{equation}
    \label{eq:r.sub.geod}
\int_0^1 g_{f(t)}^\F(\dot f(t), \dot f(t)) dt
\end{equation}
subject to the constraints $\pi(f(0)) = b_0$ and $\pi(f(1)) = b_1$. One can also show that any optimal $f(\cdot)$ must be horizontal (i.e., $\dot f(t) \in H_{f(t)}$ for all $t$). Moreover, because of the isometry between horizontal spaces, one can simplify the first constraint by choosing $f(0) = f_0$ arbitrarily in $ \pi^{-1}(b_0)$, while $f(1)$ still needs to be optimized over $\pi^{-1}(b_1)$. This remark will be used in our algorithms, because optimization over trajectories in $\F$ will generally lead to simpler formulations. 
%\todonicolas{I would be tempted to leave the following paragraph for Section 3 or 4 and not yet talk about relaxing the end time constraint in this section. Or we could also move all the rest of this section to Section 4 (?)}
%\todolaurent{I think introducing the idea here is fine, and we can be redundant and repeat it in the other sections.}
Moreover, we will often relax the constraint at $t=1$ and replace it with a penalty, minimizing
\begin{equation}
    \label{eq:r.sub.geod.2}
\int_0^1 g_{f(t)}^\F(\dot f(t), \dot f(t)) dt + U(f(1), f_1)
\end{equation}
where $f_1$ is a fixed element in $\pi^{-1}(f_1)$ and $U$ is a non-negative function blind to ``fiber'' variations, such that $U(f,f') = 0$ if and only if $\pi(f)=\pi(f')$. The definition of such functions when comparing parametrized curves or surfaces will be addressed in Section \ref{sec:chordal}.  \\

We will, in particular, consider two specific situations, both associated with group actions, for which we introduce some notation. Let $\G$ be a Lie group, i.e., both a group and a manifold, in which  the group operation is smooth. The Lie algebra of $\G$, denoted $\mf g$, is the tangent space to $\G$ at the identity. 
The left or right action of a Lie group $\G$ on a manifold $\M$ will generally be denoted with a dot: $(g, m) \mapsto g\cdot m$ (left action) and $(g,m) \mapsto m\cdot g$ (right action). This action is also assumed to be smooth as a mapping from $\G\times \M$ to $\M$. Assuming a left action, and fixing $m$, we let $A_m: g \mapsto g\cdot m$ and one defines the infinitesimal action of $\mf g$ on $\M$ by $v\cdot m = dA_m(\id)v$, $v\in \mf g$. (A similar definition can be made for right actions, with the infinitesimal action denoted $m\cdot v$.) We will also use the notation $\xi_m = dA_m(\id) $, which is a linear mapping between $\mf g$ and $T_m\M$. 
%\todonicolas{Should we reverse the order of (1) and (2) to be consistent with the order of the later sections?}

We here describe three construction involving Riemannian submersions that have been introduced for metric shape registration.  The first two of them will serve as blueprints for the constructions of Sections \ref{sec:intrinsic_metrics} to \ref{sec:hybrid}. 
\begin{enumerate}[wide]
    
    \item Assume that $\G$ is a group acting on $\F$ (on the right) and $\B$ is the quotient space 
    \[
    \F/\G = \{[f] = f\cdot \G: f\in \F\}.
    \]
%    \todomartin{Note, that this is for most of our examples not true. For non-smooth function spaces the manifold structure of the quotient space is unknown (probably wrong) and even for smooth functions $\Imm/\Diff$ is only an orbifold, which is however enough for Riemannian submersions. Maybe again to technical.}
    Assume that $\F$ and $\G$ are such that this quotient space is also a manifold. Sufficient conditions for this to hold (see, e.g., \cite{marsden2013introduction}, Chapter 9) typically require that $\G$ acts freely on $\F$, i.e., that, for all $f\in\F$, the identity $f\cdot g = f$ only holds when $g=\id$, in which case $\pi: f \mapsto f\cdot \G$ is a submersion. Assume that the metric on $\F$ is such that the action $\theta_g: f\mapsto f\cdot g$ is an isometry for all $g\in \G$. Then one can build a metric on $\B$ such that $\pi$ is a Riemannian submersion. Indeed, $\pi(f) = \pi(f')$ implies that $f' = T_g f$ for some $g$, so that $T_f\F$ and $T_{f'}\F = d\theta_g(f)T_f\F $ are isometric, and so are $V_f$ and $V_{f'}$, because $\pi\circ \theta_g = \pi$ implies that $d\pi(f') d\theta_g(f) = d\pi(f)$. This implies that $d\theta_g(f)$ also is an isometry between $H_f$ and $H_{f'}$, and that it is equal to $\psi_{f,f'}$. 
    
    Under this assumption, \eqref{eq:r.sub.geod} must be minimized subject to $f(0) = f_0$ and $f(1) = f_1 \cdot g$ for some $g\in\G$ (so that $g$ can be considered as an additional variable in the optimization). If one uses \eqref{eq:r.sub.geod.2}, the function $U$ must be such that $U(f,f') = 0$ if and only if $f' = f\cdot g$.\\

    \item For the second example, assume that $\F$ is a Lie group (with Lie algebra $\mf f$) acting transitively on the left on $\B$, so that, for all $b,b'\in \B$, there exists $f\in \F$ such that $f\cdot b=b'$. Fix an element $b_0\in B$ and define $\pi:\F\to \B$ by $\pi(f) =f\cdot  b_0$, which is a submersion. If one assumes that the metric on $\F$ is right-invariant, so that transformations $\theta_f: f'\mapsto f'f$ are isometries, then horizontal spaces are isometric. Indeed, the left-invariance assumption ensures that the tangent spaces to $\F$ are isometric, with $T_f\F = d\theta_f(\id) \mathfrak f$. Moreover, it is easy to check that vertical spaces are right translations of the spaces 
    \[
    \mathfrak v_b = \{h\in \mathfrak f: h\cdot b = 0\}
    \]
%    \todonicolas{I wonder if the notations are not going to end up being confusing if $q$ denotes here a diffeomorphism of $\R^n$. I would rather write $g$ or $\varphi$ instead and keep $q$ for a parametrized immersion which should be the state of the optimal control problems in both formulations after all the reductions.}
   % \todolaurent{OK, using $Q$ for the top space was a bad idea. Any suggestion for another choice? Maybe $\Theta$ or $\Gamma$?}
    %\todomartin{How about $\mathcal N$? Or even just $N$}
    through $V_f = d\theta_f(\id) \mathfrak v_b$ when $\pi(f) = b$. This in turn implies that  $H_f = d\theta_f(\id) \mathfrak h_b$ with $\mf h_b = \mathfrak v_b^\perp$, proving that horizontal spaces are all isometric. Moreover, if $\pi(f) = \pi(f')$, then $\pi\circ \theta_f = \pi\circ \theta_{f'}$ and $d\pi(f)d\theta_f(\id) = d\pi(q') d\theta_{f'}(\id)$ so that $\psi_{f,f'} d\theta_f(\id) = d\theta_{f'}(\id)$ and $\psi_{f,f'}$ is an isometry.\\
    
    The right-invariance property of the group metric provides an alternate formulation of \eqref{eq:r.sub.geod} or \eqref{eq:r.sub.geod.2}. Indeed, this property implies that
    \[
    g_f^\F(\dot f, \dot f) = g_{\id}^\F(d\theta_{f^{-1}}(f)\dot f, d\theta_{f^{-1}}(f)\dot f)
    \] and introducing the new variable $v(t) = d\theta_{f(t)^{-1}}(f)(\dot f(t))\in \mf f$, \eqref{eq:r.sub.geod} can be expressed in the form of an optimal control problem, namely minimizing
    \begin{equation}
    \label{eq:r.sub.oc}
\int_0^1 g_{\id}^\F(v(t), v(t)) dt
\end{equation}
subject to $f(0) = f_0$, $\dot f = \theta_f v$ and $\pi(f(1)) = b_1$ (where $f_0$ is chosen so that $\pi(f_0)=b_0$). A similar transformation can be done with \eqref{eq:r.sub.geod.2}, in which it suffices to take the function $U$ in the form $U(f,f') = \tilde U(f\cdot b_0, f'\cdot b_0)$, where $\tilde U$ is defined on $\B\times \B$,  to obtain the required invariance. This leads to the formulation that is used 
 in Section \ref{sec:LDDMM} for the introduction of outer deformation metrics.\\
 
 \item In our third example, we assume that $\G$ is a Lie group with a left action on $\B$, and take $\F = \G \times \B$. (We do not require the action to be transitive.) We consider the mapping $\pi: \F\to\B$ such that $\pi(g,\beta) = g\cdot \beta$, which is a submersion. If one equips $\F$ with a Riemannian metric that is invariant by the right action of $\G$ on $\F$ defined by $(h,\beta)\cdot g = (hg, g^{-1}\cdot\beta)$, which leaves fibers invariant since $\pi((h,\beta)\cdot g) = \pi(h,\beta)$, then it is easy to show that horizontal spaces over a given point $b\in \B$ are isometric, therefore providing a new Riemannian metric. This construction corresponds to the metamorphosis metric introduced in \cite{miller2001group,trouve2005metamorphoses}, which can be expressed as
 \[
 g_b^\B(\eta, \eta) = \min\{G_{(\id, b)}^\F((v, h), (v, h)): \eta = v\cdot b + h \}\,.
 \]\\
\end{enumerate}
%\todonicolas{Should we maybe reverse the order of (1) and (2) since we then start by discussing intrinsic metrics before outer metrics?}

Notice that one can use (1) after (2) in sequence where $\F$ acts transitively on $\B$ and $\G$ has another (free) action on $\B$ also. It is then necessary that the norm obtained after (1) be invariant by $\G$ in order for (2) to apply. We will review examples in Section \ref{sec:LDDMM}. 

% Metrics constructed through Riemannian submersions:
% $$ \pi: \Diff(\R^n) \rightarrow \mathcal{S}$$
% with $\mathcal{S}=\Diff(\R^n).q_0$ (outer metrics). Alternatively,
% $$ \pi: \Imm(M,\R^n) \rightarrow \Imm(M,\R^n)/\Diff(M) $$
% leads to inner metrics. 

% Possible metamorphosis model through the submersion
% $$ \pi: \Diff(\R^n) \times \Imm(M,\R^n) \rightarrow \Imm(M,\R^n)/\Diff(M) $$
% $$ (\phi,q) \mapsto (\phi \circ q)\circ \Diff(M)$$

%\section{Measuring discrepancy between curves and surfaces}
%Hausdorff distance, current/varifold...

%can "crudely" quantify shape distance but with no underlying dynamical evolution model.

\section{Optimal control framework}\label{sec:control}
As we will evidence in the sections to follow, there is a particular advantage in interpreting the variational problems appearing in metric shape registration as optimal control problems, which gives a natural and convenient framework to derive optimality equations or design algorithms to numerically solve them. In this section, we give a brief overview of the main principles of optimal control, which will be then applied to the cases of intrinsic, outer and hybrid metric matching in Sections \ref{sec:intrinsic_metrics}, \ref{sec:LDDMM} and \ref{sec:hybrid}. We will however restrict our presentation to systems governed by finite-dimensional ODEs. We should yet point out that, unless discretized, the optimal control problems appearing in shape analysis are generally infinite dimensional, which adds specific difficulties.  We refer readers to \cite{Arguillere15} for a partial extension of the results mentioned below to the infinite dimensional case.       

We assume that the state of the system at any time $t$ is modeled by a variable $x(t)\in \R^n$ that evolves according to an ODE of the form $\dot{x}(t) = f(x(t),u(t))$ with $u(t) \in \Omega \subset \R^k$ being the control variable, $\Omega$ its constraint set, and $f: \R^n \times \Omega \rightarrow \R^n$ a function satisfying the adequate regularity assumptions for solutions of the ODE to exist. The class of optimal control problem we consider consists in determining a control $t \mapsto u(t)$ that brings the system from its initial state $x(0)=x_0$ to a final state $x(1)=x_1$ with $x_1$ belonging to a certain target set $\mathcal{C}_1$ (e.g a single point or a submanifold of $\R^n$) while minimizing a cost functional of the form:
\begin{equation}
\label{eq:opt_contr_cost}
    C(u) = \int_0^1 L(x(t),u(t)) dt + U(x(1))
\end{equation}
in which $L: \R^n \times \Omega$ is called the Lagrangian functional of the system while $U(x(1))$ is a cost term on the end time state. The two typical situations we encounter in this chapter are the cases where $\mathcal{C}_1=\{x_1\}$ with $U=0$ which corresponds to a boundary value (or Lagrange) problem or $\mathcal{C}_1 = \R^n$ with $U(x(1))$ measuring the distance of $x(1)$ to the target set which corresponds to a relaxed (or Mayer-Lagrange) problem.\\

There are several questions that optimal control theory attempts to address, the first of which being the existence of solutions. This existence problem usually involves first the existence of admissible trajectories between the initial state $x_0$ and the target set $\mathcal{C}_1$: this is the issue of \textit{controllability} and is particularly important in the case of boundary value problems. It is however a difficult question to address in general and often only local controllability can be established. Conditions of global controllability exist nonetheless in the cases of linear or sub-Riemannian control systems (c.f.,  for instance, the presentation in \cite{Coron2007}). Next, one needs to show that the minimum of $C$ over all admissible controls is achieved. This typically relies on a compactness and lower semi-continuity argument on the cost functional, for which some conditions have to be carefully checked with each specific problem at hand.

The second point of optimal control is to provide some characterization of the solutions. To that end, a very powerful result is the \textbf{Pontryagin maximum principle} (PMP), first introduced in \cite{Pontryagin1962}. For the control system and cost functional considered here, it can be written as follows. We first define the \textit{Hamiltonian} $H_{\lambda}: \R^n \times \R^n \times \Omega \rightarrow \R$ as:
\begin{equation}
    \label{eq:opt_control_Hamil}
    H_{\lambda}(x,p,u) = p^T f(x,u) - \lambda L(x,u)
\end{equation}
where $p$ denotes the co-state variable that we will often call \textit{momentum}. The PMP states that if $u_*$ is an optimal control then there must exist $\lambda \in \{0,1\}$ and a time-dependent momentum $p:[0,1] \rightarrow \R^n$ such that $(\lambda,p)\neq(0,0)$ and that the following Hamiltonian system of equations is satisfied:
\begin{equation}
    \label{eq:opt_control_Hamil_eq}
\left\lbrace\begin{aligned}
&\dot{x}(t)=\partial_p H_\lambda(x(t),p(t),u_*(t)),\\
&\dot{p}(t)=-\partial_x H_\lambda(x(t),p(t),u_*(t)) \\
&u_*(t)=\mathrm{argmax}_{u\in\Omega}H_{\lambda}(x(t),p(t),u).
\end{aligned}\right.    
\end{equation}
Moreover, the following transversality condition holds: $p(1) + \lambda \nabla U(x(1))$ is orthogonal to the tangent space $T_{x(1)}\mathcal{C}_1$ of the constraint set at $x(1)$.

There are thus two possible types of solutions for such optimal control problems. When $\lambda =0$, $u_*$ and the optimal trajectory are independent of the Lagrangian: these are called \textit{abnormal solutions} of the system. Note that such solutions are rather singular and exist only for special systems and/or choice of boundary conditions, and we generally not consider such solutions in practice. For $\lambda =1$, we get \textit{normal solutions}. The transversality condition is directly related to the target set $\mathcal{C}_1$: in particular, in the unconstrained case $\mathcal{C}_1 = \R^n$, it reduces to $p(1) = -\nabla U(x(1))$ whereas in the boundary problem $\mathcal{C}_1=\{x_1\}$, we get no conditions on $p(1)$ (but instead $x(1)=x_1$).\\

The third equation in \eqref{eq:opt_control_Hamil_eq} may, in some cases, be used to express the optimal control in feedback form, namely as a function of the state and the co-state, i.e as $u_*(t)=\eta(x(t),p(t))$. For instance, if $\Omega$ is an open subset, one can write first order optimality equations $\partial_u H_{\lambda}(x(t),p(t),u_*(t)) = 0$ and attempt to explicit $u_*$ from the implicit equations. In this case, the Hamiltonian system can be rewritten in more compact form. Indeed, if we denote by $\tilde{H}(x,p) = p^T f(x,\eta(x,p)) - L(x,\eta(x,p))$ the \textit{reduced Hamiltonian}, \eqref{eq:opt_control_Hamil_eq} becomes:
\begin{equation}
    \label{eq:opt_control_reduced_Hamil_eq}
\left\lbrace\begin{aligned}
&\dot{x}(t)=\partial_p \tilde{H}(x(t),p(t)),\\
&\dot{p}(t)=-\partial_x \tilde{H}(x(t),p(t)) 
\end{aligned}\right.    
\end{equation}
These reduced Hamiltonian equations eventually show that (normal) optimal solutions are fully characterized by the initial momentum $p_0 \in \R^n$. This leads to a considerable reduction of the search space which can be exploited in numerical approaches.  \\

In fact, several different algorithmic strategies can be implemented to approximate solutions of the previous problems. We briefly discuss two among the main ones. A first approach, commonly known as \textbf{trajectory optimization}, is to start from a discretization in time of all the problem variables. Considering time samples $0=t^{(0)}<t^{(1)}<\ldots<t^{(N)}=1$ and the associated collocation points for the state $x^{(0)}, x^{(1)},\ldots,x^{(N)}$ and control $u^{(0)},u^{(1)},\ldots,u^{(N)}$, the principle is to define a numerical scheme (such as Euler, Runge-Kutta,\ldots) for the ODE $\dot{x}(t)=f(x(t),u(t))$ and a quadrature rule to approximate the cost functional \eqref{eq:opt_contr_cost} in order to rewrite both with respect to the $x^{(i)}$ and $u^{(i)}$. This turns the optimal control problem into a nonlinear programming one, which can be tackled with standard constrained optimization solvers such as interior point method or sequential quadratic programming. Trajectory optimization methods have the advantage of being easy to design and relatively robust to the initialization of the state and control variables. However, in order to get a good accuracy for the solution, they can often result in very high-dimensional optimization problems. 

A second possible approach is to exploit, when available, the dimensionality reduction given by the PMP and in particular the reduced Hamiltonian equations \eqref{eq:opt_control_reduced_Hamil_eq} by searching directly for an optimal initial momentum $p_0$. This is the principle of \textbf{shooting} methods. To keep this section short, we will only say that those basically amount in optimizing the cost functional with respect to $p_0$ with first or second order descent schemes, the main difficulty being to compute the gradient of \eqref{eq:opt_contr_cost} as a function of $p_0$. Numerical or automatic differentiation techniques can be very often used for that purpose. Another common strategy is to obtain the gradient by numerically solving the adjoint Hamiltonian system of equations. The optimization problem to solve is typically of much lower dimension than with trajectory optimization, although shooting methods may be more sensitive to the initialization of $p_0$ and the presence of local minima.

\section{Chordal Metrics on Shapes}
\label{sec:chordal}
\subsection{Motivation}
The Riemannian metric framework introduced in Section \ref{sec:submersion} leads to geodesic distances between two given shapes $b_0$ and $b_1$ that are formally expressed through a minimization over paths connecting $q_0 \in \pi^{-1}(b_0)$ to $q(1) \in \pi^{-1}(b_1)$. As already mentioned, it is often necessary to relax the end-time constraint for multiple reasons. 

Indeed, in certain situations, there may not exist paths connecting $q_0$ to an element of the fiber $\pi^{-1}(b_1)$. Moreover, one could argue more generally that many variations between shapes may be the result of undesirable perturbations like acquisition noise, segmentation issues or boundary effects. Imposing exact matching constraint is then very likely to result in distances with poor robustness properties to such perturbations and thus of little interest in problems such as  statistical shape classification on real datasets. 

A second issue is the complexity of optimizing \eqref{eq:r.sub.geod} under the constraint that $q(1)$ belongs to $\pi^{-1}(b_1)$, which requires an additional exploration of vertical motions in the fiber of $b_1$. Those consist essentially in reparametrizations of the shape $b_1$. While numerical frameworks for discretizing and optimizing over $\Diff(M)$ have been quite commonly used for closed and open curves (i.e. when $M=\S^1$ or $M=[0,1]$) in works like \cite{Klassen2004,Jermyn2011,BBHM2017,Srivastava2016}, these are harder to extend to cases like closed or open surfaces and with the general class of metrics we are considering here. In this direction, only work for a particular metric of order one on the space of sphere-like surfaces has been done in \cite{kurtek2013landmark,laga2017numerical}.

In this section, we thus seek an auxiliary family of shape distances to serve as a relaxation of the constraint $q(1) \in \pi^{-1}(b_1)$, which would be 1) defined explicity, 2) reparametrization-blind by construction and 3) easy/fast to compute numerically. These shape distances are based on the idea of "embedding" shape spaces into certain functional spaces, in contrast with the previous submersion setting, and we thus refer to them as chordal metrics. %Although it may seem that both chordal and submersion distances attempt to eventually address the same problem, these two approaches should be thought as complementary rather than redundant with one another, as we shall highlight in the rest of this section.

\subsection{General Principle}
Chordal metrics are constructed by mapping $\Imm(M,\R^n)$ into a certain Banach space $(\mathcal{H},\|\cdot\|_{\mathcal{H}})$ (typically a space of distributions or measures over a certain feature set as we detail below). The norm on $\mathcal{H}$ can be then pulled back to $\Imm(M,\R^d)$ based on this mapping. 

\begin{figure}
\centering
 \includegraphics[width=10cm]{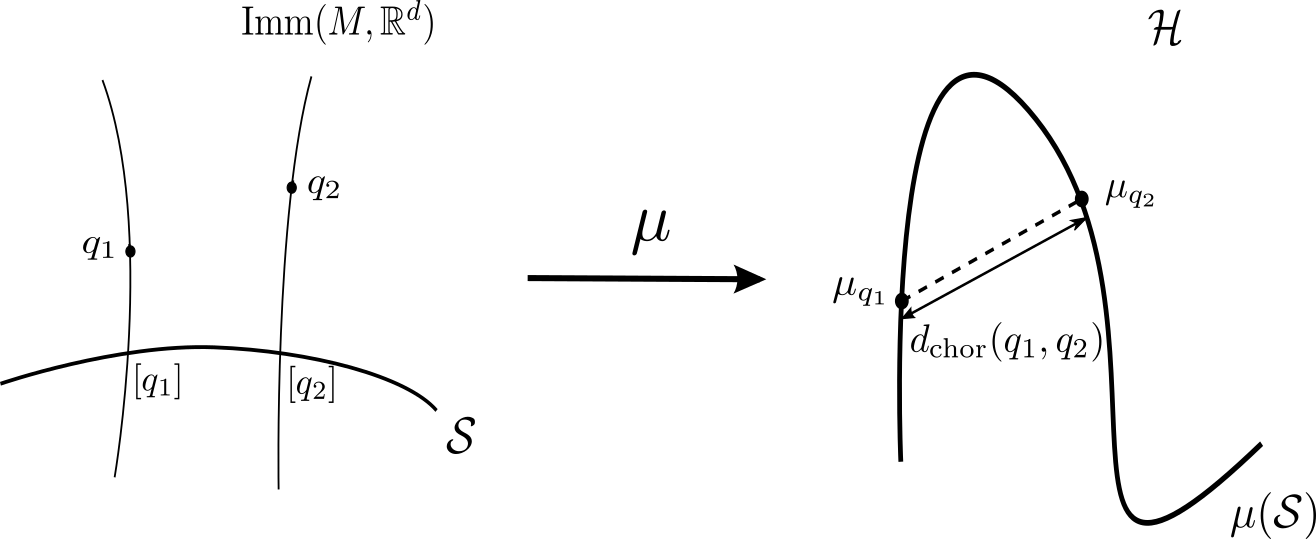}
 \caption{Shape space embedding and chordal distance.}\label{fig:chordal_dist}
\end{figure}

More specifically let $\mu: \Imm(M,\R^d) \rightarrow \mathcal{H}$ be a mapping such that the representation of any $q \in \Imm(M,\R^d)$ as $\mu_{q} \in \mathcal{H}$ is \textit{independent of parametrization}, that is for all $\tau \in \Diff(M)$ we have $\mu_{q \circ \tau} = \mu_q$. Then, one has a well-defined quotient mapping $[\mu]$ on $\mathcal{S}=\Imm(M,\R^d)/\Diff(M)$ and we can introduce for $q_1,q_2 \in \Imm(M,\R^d)$:
\begin{equation}
\label{eq:chordal_dist}
    d_{\text{chor}}([q_1],[q_2]) = d_{\text{chor}}(q_1,q_2) = \|\mu_{q_1} - \mu_{q_2}\|_{\mathcal{H}}
\end{equation}
which is precisely the chordal distance in $\mathcal{H}$ on the image of $\Imm(M,\R^d)$ by $\mu$, as illustrated in Figure \ref{fig:chordal_dist}. However, \eqref{eq:chordal_dist} is in general only a pseudo-distance on $\mathcal{S}$. To obtain a true distance, one needs the additional assumption that $[\mu]$ is injective, in other words that the space of shapes is embedded in $\mathcal{H}$. 

Note that the resulting (pseudo-)distances are defined between any two elements of $\mathcal{S}$ but are fundamentally different from the Riemannian or sub-Riemannian metrics of Section \ref{sec:submersion} since, for chordal distances, there are no corresponding geodesics on the shape space. On the other hand, with adequate choices of $\mu$ and of the norm $\|\cdot\|_{\mathcal{H}}$, such distances can still provide a reasonable proximity criterion which can be computed in a much more direct way.         

\subsection{Oriented varifold distances}
There are multiple possible constructions of a mapping $\mu$ satisfying the conditions above. A popular approach, inspired from the field of geometric measure theory, is to represent shapes as measures or distributions and compare them with the corresponding dual metrics. This is the basic idea behind such frameworks like measures \cite{Glaunes2004}, currents \cite{Glaunes2008}, varifolds \cite{Charon2013}, and recently the works of \cite{Roussillon2016} on normal cycles and \cite{Feydy2017} on optimal transport discrepancies. In the rest of this section, we focus on the setting of \textit{oriented varifold} for curves and surfaces presented in \cite{Charon2017}, which has the advantage of encompassing several of these distance families in one single model. We will also restrict the presentation to the case where shapes are either curves or hypersurfaces in $\R^d$, i.e $\text{dim}(M)=1$ or $\text{dim}(M)=d-1$. 

We let $\mathcal{H}=W^*$ with $W$ being a Hilbert space of test functions on $\R^d \times \S^{d-1}$, with the continuous embedding $W \hookrightarrow C_0(\R^d \times \S^{d-1})$, and define $\mu$ as follows. For any $q \in \Imm(M,\R^d)$, $\mu_q$ is the distribution in $W^*$ such that for all $\omega \in W$:
\begin{equation}
    \label{eq:def_mu_orvar}
    \mu_q(\omega) = \int_{M} \omega(q(m),\vec{t}(m)) d\vol_q(m) 
\end{equation}
where $\vec{t}(m) \in \S^{d-1}$ denotes the unit oriented tangent (resp. normal) vector to the curve (resp. hypersurface) at $q(m)$, and $\vol_q$ is the intrinsic volume density of $q$ (in the case of curves, this would correspond to the integration with respect to arc length). A simple change of variable shows that $\mu_q$ is independent of the parametrization. Intuitively, \eqref{eq:def_mu_orvar} amounts in representing a shape as a distribution of its tangent spaces attached to the different positions of its points.

We can then introduce the associated chordal pseudo-distance in the dual Hilbert space $W^*$ which we shall write:
\begin{equation}
\label{eq:chordal_dist_orvar}
    d_{\text{var}}([q_1],[q_2])^2 = \|\mu_{q_1} - \mu_{q_2}\|_{W^*}^2 = \|\mu_{q_1}\|_{W^*}^2 -2\langle \mu_{q_1},\mu_{q_2} \rangle_{W^*} + \|\mu_{q_2}\|_{W^*}^2.
\end{equation}
Such a distance depends obviously of the choice of the test space $W$ and its metric. Yet the main interest of this model is that $d_{\text{var}}$ can be in fact expressed explicitly in the right setting. Indeed, given the assumption $W\hookrightarrow C_0(\R^d \times \S^{d-1})$, $W$ is a reproducing kernel Hilbert space (RKHS) on $\R^d \times \S^{d-1}$ which implies that $W$ has a unique reproducing kernel $K_W: (\R^d \times \S^{d-1})^2 \rightarrow \R $. The reproducing kernel property allows one to write $\|\cdot\|_{W^*}$ in terms of $K_W$ and one can show in particular that 
\begin{equation*}
    \|\mu_{q} \|_{W^*}^2 = \iint_{M\times M} K_W((q(m),\vec{t}(m)),(q(m'),\vec{t}(m'))) d\vol_q(m) d\vol_q(m').
\end{equation*}

In practice, we proceed in the reverse way and start from a positive kernel on $\R^d \times \S^{d-1}$ that implicitly specifies a space $W$. We generally restrict to separable kernels, such that $K_W$ is the product of a positive kernel on $\R^d$ and of a positive kernel on $\S^{d-1}$. Furthermore, it is often relevant to seek metrics that are invariant to the action of rotations and translations. This can be achieved with kernels of the following form:
\begin{equation*}
    K_W((x,\vec{t}),(x',\vec{t}')) = \rho(|x-x'|^2) \gamma(\scp{\vec{t}}{\vec{t}'}),
\end{equation*}
which are functions only of the distance between positions $x$ and $x'$ and of the angle between oriented unit vectors $\vec{t}$ and $\vec{t}'$, with $\rho: \R_+ \rightarrow \R$ and $\gamma: [-1,1]\rightarrow \R$. Given two such functions $\rho$ and $\gamma$ that define positive kernels on $\R^d$ and $\S^{d-1}$ respectively, the oriented varifold metric becomes fully explicit:
\begin{equation}
\label{eq:oriented_var_kernel_norm}
    \|\mu_{q} \|_{W^*}^2 = \iint_{M\times M} \rho(|q(m)-q(m')|^2) \gamma(\scp{\vec{t}(m)}{\vec{t}(m')}) d\vol_q(m) d\vol_q(m').
\end{equation}

As mentioned earlier, $d_{\text{Var}}$ defined by \eqref{eq:chordal_dist_orvar} and \eqref{eq:oriented_var_kernel_norm} is only a pseudo-distance on $\mathcal{S}$ unless the RKHS $W$ associated with the kernel is such that the mapping $[\mu]$ is injective. Unfortunately, it can be seen that, in the case of oriented varifold representation \eqref{eq:def_mu_orvar}, this is actually impossible for the full space of unparametrized immersed curves, no matter the choice of $W$ (see for example \cite{BauerCoCV2018} for a counter-example with curves). However, sufficient conditions exist for $d_{\text{var}}$ to be a distance on the subset of \textit{embedded} submanifolds $\Emb(M,\R^d)/\Diff(M)$.
\begin{theorem}
\label{thm:distance_kernel_var}
Assume that $\rho$ and $\gamma$ are $C^1$-functions, that the kernel defined by $\rho$ is $C_0$-universal on $\R^d$ and that $\gamma(1)>0$, $\gamma(-u) \neq \gamma(u)$ for all $u\in[-1,1]$. Then $d_{var}$, defined by \eqref{eq:chordal_dist_orvar} and \eqref{eq:oriented_var_kernel_norm}, is a distance on the set of embedded shapes.
\end{theorem}
The $C_0$-universality assumption is essentially a density property of the RKHS associated with $\rho$ (cf \cite{Carmeli2010} for details). It is satisfied by families of kernels like Gaussian $e^{-\frac{|x-x'|^2}{\sigma^2}}$, Cauchy $1/\left(1+\frac{|x-x'|^2}{\sigma^2}\right)$, Wendland... The assumptions on $\gamma$ are much milder: typical choices include for instance the linear kernel $\gamma(\langle \vec{t},\vec{t}' \rangle) = \langle \vec{t},\vec{t}' \rangle$ or the restricted Gaussian kernel on $\S^{d-1}$ $\gamma(\langle \vec{t},\vec{t}' \rangle) = e^{\frac{2}{\sigma^2}(1-\langle \vec{t},\vec{t}' \rangle)}$. 

Although we do not obtain, with this framework, a true distance on $\mathcal{S}$, we point out that the result of Theorem \ref{thm:distance_kernel_var} may in fact still hold with the same assumptions for a larger class than embedded shapes: for instance, it was extended to immersed curves with transverse self-intersections in \cite{BauerCoCV2018}. In all cases, and for the purpose of this work, $d_{\text{var}}$ gives in practice a simple and explicit notion of shape proximity to be used as a relaxation of the exact matching constraint.

\subsection{Numerical aspects}
In practical cases, shapes are discretized and oriented varifold metrics are computed based on approximations of the integrals appearing in \eqref{eq:oriented_var_kernel_norm}. Although there exist different possible discretization schemes for curves and surfaces, a common situation is to deal with vertices and a mesh structure such as lists of segments or triangles. In such cases, we can view the parameter space $M$ as a reunion of simplices $M = \bigcup_{i=1}^{N_S} M_{i}$ where $M_i$ are disjoint faces delimited by vertices $m_{i,1},\ldots,m_{i,p}$ ($p=2$ for segments, $p=3$ for triangles). Then \eqref{eq:oriented_var_kernel_norm} may be rewritten as:
\begin{equation*}
    \|\mu_{q} \|_{W^*}^2 = \sum_{i,j=1}^{N_S} \iint_{M_i\times M_j} \rho(|q(m)-q(m')|^2) \gamma(\scp{\vec{t}(m)}{\vec{t}(m')}) d\vol_q(m) d\vol_q(m')
\end{equation*}
and we can approximate every integral on the product simplex $M_i \times M_j$ using a certain quadrature rule. For example, the approach proposed in \cite{Charon2013,Charon2017} consists in simply approximating the integrals by their interpolated value at the center of simplices. Specifically, denoting $\bar{q}_i$ the barycenter of the vertices $(q(m_{i,1}),\ldots,q(m_{i,p}))$, $\vec{t}_i$ the tangent or normal orientation vector to this simplex and $\vol_q(M_i)$ the length or area of the simplex $q(M_i)$, the previous norm is approximated as:
\begin{equation}
\label{eq:oriented_var_kernel_norm_approx}
    \|\mu_{q} \|_{W^*}^2 \approx \sum_{i,j=1}^{N_S} \rho(|\bar{q}_i-\bar{q}_j|^2) \gamma(\scp{\vec{t}_i}{\vec{t}_j}) \vol_q(M_i) \vol_q(M_j).
\end{equation}
It can be shown that the error in \eqref{eq:oriented_var_kernel_norm_approx} is controlled by the maximum diameter of the simplices. In addition, under certain assumptions, convergence of discrete to continuous varifold shape representations has been established. From a numerical standpoint, the computational cost of \eqref{eq:oriented_var_kernel_norm_approx} or of its gradient with respect to vertex positions is typically quadratic in the number of simplices. It is however highly parallelizable and many recent implementations in CUDA take advantage of GPU architectures. We refer the interested reader to \cite{ArguillereSIIMS16,Charon2017,CCGRG2018} for more detailed discussions on discrete varifold approximations and computations.

\section{Intrinsic metrics}
\label{sec:intrinsic_metrics}
In this section we will study reparametrization-invariant Sobolev metrics on the space of parametrized shapes 
$\Imm(M,\R^d)$. The reparametrization invariance of the metric will allow us to consider induced metrics on the space of unparametrized shapes as described in Section~\ref{sec:submersion} paragraph (1).
We call these metrics \emph{intrinsic} as they will be defined solely in terms of 
the shape (immersion, resp.), which is in  contrast to the \emph{outer deformation metric models} studied in Section~\ref{sec:LDDMM}, that are defined using the geometry of the ambient space of the (parametrized) shapes. In order to simplify the presentation we will first introduce the class of intrinsic metrics for the space of infinitely smooth immersions and will later show that they can be extended to smooth Riemannian metrics on the bigger space of all immersions of (sufficiently high) Sobolev class.

\subsection{Reparametrization-invariant metrics on parametrized shapes}
To view $\Imm_{C^{\infty}}(M,\R^d)$ as an (infinite-dimensional) manifold we observe that it is an open subset of the vector space  of all smooth functions $C^{\infty}(M,\R^d)$. Thus  tangent vectors to an immersion $q$ -- vector fields along the immersion-- can be identified with  smooth function $h\in C^{\infty}(M,\R^d)$, see \cite{Michor1991}. To define a Riemannian metric on this infinite dimensional manifold we have to specify a family of inner products
\begin{equation}
\Imm_{C^{\infty}}(M,\R^d)\times C^{\infty}(M,\R^d)\times C^{\infty}(M,\R^d)\ni (q,h,k) \mapsto G_q(h,k)\in \mathbb R\;,    
\end{equation}
that depends smoothly on the footpoint $q\in \Imm_{C^{\infty}}(M,\R^d)$. As our main goal is the study of the space of unparametrized shapes we will always require the metric to be invariant under the action of the group of reparametrizations, i.e.,
\begin{equation}
G_q(h,k)=G_{q\circ \tau}(h\circ\tau,k\circ\tau)\;,    
\end{equation}
for all $\tau \in \Diff_{C^{\infty}}(M)$. This assumption is  a necessary condition for a metric $G$ to induce a Riemannian metric on the quotient space $\mathcal S$ such that the projection $\pi$ is a Riemannian submersion. In this infinite-dimensional situation one also has to prove the existence of the horizontal space, as defined in Section~\ref{sec:submersion}. For all the metrics studied in this section this second condition will be satisfied, which is in contrast to  metrics on $\Imm_{C^{\infty}}(M,\R^d)$, that are induced via outer deformation metric models, where the horizontal bundle only exists in some Sobolev completion of the tangent bundle~\cite{micheli2013sobolev}. 

We are now ready to define the specific class of metrics, that we will study in this section. Therefore we let $\Delta_q$ be the Laplace operator\footnote{To formally define the Laplace operator let $g=q^*\langle.,.\rangle$ be the induced surface metric with corresponding covariant derivative $\nabla$. Then the induced Laplacian of a function
$h\in C^{\infty}(M,\mathbb R^d)$ is given by $\Delta_q h=\operatorname{tr}(g^{-1}\nabla^2 h)$. It will be important for us, that the Laplacian is equivariant under the action of the diffeomorphism group, i.e., $\Delta_{q\circ \tau}(h\circ \tau)=(\Delta_q h)\circ\tau$ for all $\tau \in \Diff_{C^{\infty}}(M)$.} induced by the submanifold $q$. A reparametrization-invariant Sobolev metric of order $s$ is given by:
\begin{equation}\label{eq:RiemMetric}
G_q(h,k)=\int_M \langle (1+\Delta_q)^{s/2} h, (1+\Delta_q)^{s/2} k\rangle_{\mathbb R^d}\operatorname{vol}_q\;,
\end{equation}
where $\operatorname{vol}_q$ denotes integration with respect to the intrinsic volume density of the submanifold $q$.
For odd (or non-integer) $s$ we define the fractional Laplacian using  analytic
functional  calculus, see e.g. \cite{bauer2018smoothness}. 
For closed manifolds $M$ we can use integration by parts to rewrite the formula of the metric as
\begin{equation}
G^s_q(h,k)=\int_M \langle (1+\Delta_q)^{s} h, k\rangle_{\mathbb R^d}\operatorname{vol}_q\;.
\end{equation}
This class of metrics has been first introduced in the case of planar curves in \cite{Michor2007,Sundaramoorthi2011,Mennucci2008,Sundaramoorthi2008} and has been later generalized to immersions between general Riemannian manifolds in \cite{Bauer2011b,bauer2012sobolev,Bauer2014}. 

The invariance of the metric under the action of the reparametrization group follows directly from the equivariance of the Laplacian and the transformation formula for multidimensional integration. 
\begin{remark}\label{rem:invariances}
It is worth to note that the metric is in addition invariant under the Euclidean motion group (translations and rotations) and thus induces a metric on the space of shapes modulo Euclidean motions. It is however not invariant under the action of the scaling group. To achieve a scale invariant version one can add weights depending on the total volume of the surface into the definition of the metric, see \cite{bauer2012sobolev,Michor2008a,Mennucci2008}. 
\end{remark}

More generally, one can consider metrics defined in terms of  fields of (pseudo)-differential operators $L_q$:
\begin{equation}
G^L_q(h,k)=\int_M \langle L_q h, L_q k\rangle_{\mathbb R^d}\operatorname{vol}_q\;,
\end{equation}
where $L_q$ is for each $q\in \Imm(M,\mathbb R^d)$ a positive, (pseudo) differential operator, that is equivariant under the action of the reparametrization group. The equivariance property is necessary to ensure the required invariance of the induced metric.

The following result states that Sobolev metrics as defined above can be extended to smooth Riemannian metrics on spaces of immersions of finite regularity:
\begin{theorem}
For any $s\in \mathbb N$, $\mathbb R\ni r\geq s$ and $r>\frac{d}2+1$ the Sobolev metric of order $s$ extends to a smooth Riemannian metric on the space of Sobolev immersions 
$\Imm_{H^r}(M,\mathbb R^d)$ of order $r$. 
\end{theorem}

\begin{remark}[Weak and strong Riemannian metrics]
In infinite dimensions there exists two types of Riemannian metrics, which are referred to as weak and strong metrics:
a Riemannian metric on an infinite dimensional manifold $\mathcal F$ is called strong if the induced topology of the geodesic distance function corresponds with the original manifold topology \cite{Lang1999}.
As a consequence one can view the metric as a bijective mapping $G:T\mathcal F \to T^*\mathcal F$. For weak Riemannian metrics the topology of the metric is strictly weaker than the standard topology and the mapping $G:T\mathcal F \to T^*\mathcal F$ is only injective but not surjective. As a consequence many standard results of finite-dimensional geometry (in which all metrics are strong) may fail in this situation, e.g. the geodesic spray might not exist and the induced geodesic distance can vanish identically on the manifold. We will see an example of such an ill-behavior later in Section~\ref{sec:geodesicdist} in the context of the $L^2$-metric. Sobolev metrics on spaces of immersions are usually only weak Riemannian metrics. The only exception being the metric of order $s>\frac{d}2+1$ on the space of Sobolev immersions $\Imm_{H^s}(M,\mathbb R^d)$ of the same order. 
\end{remark}
\begin{example}
To illustrate the above defined metrics, we want to present the metric of order one in the special situation $\Imm(\S^1,\mathbb R^2)$, i.e., for the case of metrics on the space of planar curves.   In that case the induced volume density of a curve $c$ is given by arc-length integration $\operatorname{vol}_c=|\partial_{\theta} c|d\theta$ and the Laplacian of the curve $c$ is given in terms of arc-length differentiation $\Delta_c=-\frac{1}{|\partial_{\theta} c|}\partial_\theta \frac{1}{|\partial_{\theta} c|}\partial_\theta $. Thus the metric can be written as:
\begin{align*}
G^1_c(h,k)&=\int_M \left\langle \left(1-\frac{1}{|\partial_{\theta} c|}\partial_\theta \frac{1}{|\partial_{\theta} c|}\partial_\theta\right) h, k\right\rangle_{\mathbb R^2} |\partial_{\theta} c|\, d\theta\\
&=\int_M \left(\left\langle h, k\right\rangle_{\mathbb R^2}|\partial_{\theta} c| +\frac{1}{|\partial_{\theta}c|} \left\langle   \partial_\theta h, \partial_{\theta} k\right\rangle_{\mathbb R^2}\right)\, d\theta\;.
\end{align*}
Using the Sobolev multiplication theorem it is easy to see that this metric
extends to the Sobolev completion $\Imm_{H^r}(M,\mathbb R^d)$ for $r>\frac{d}2+1$.
\end{example}
\subsection{The metric on the space of unparametrized shapes}
As described in Section~\ref{sec:submersion} the metric on the quotient space $\mathcal S$ corresponds to the metric on the space of immersions restricted to horizontal tangent vectors. To determine the horizontal space we first have to calculate the vertical space of the submersion $\pi: \Imm_{C^{\infty}}(M,\mathbb R^d)\to \mathcal S$: therefore we recall that elements in $\Imm_{C^{\infty}}(M,\mathbb R^d)$ that differ only by their parametrization are identified on the space of unparametrized shapes. Tangent vectors $h$ that only change the parametrization of $q$ are vector fields that are always tangent to the surface, i.e.,
\begin{equation*}
V_q=\left\{h=dq.X:  X\in C^{\infty}(M,TM)\right\}\;.    
\end{equation*}
To determine the horizontal space we have to calculate the orthogonal complement of $V$ with respect to the metric $G$. It is tempting to believe that horizontal tangent vectors are all vectors that are pointwise normal to the surface. However, it turns out that this is not true in general. Instead one has to solve a partial differential equation to calculate the horizontal path of a tangent vector:
\begin{lemma}[\cite{Bauer2011b}]
Let $M$ be a closed manifold. Then the horizontal bundle of the projection $\pi: \Imm_{C^{\infty}}(M,\mathbb R^d)\to \mathcal S$ with respect to the Sobolev metric $G$ of order $s$ is given by
\begin{equation*}
H_q=  \left\{h\in T_q \Imm_{C^{\infty}}(M,\mathbb R^d): \langle (1+\Delta)^{2s} h,dq.X\rangle=0,\; \forall X\in C^{\infty}(M,TM)\right\}\;.   
\end{equation*}
\end{lemma}
In cases where the horizontal projection of tangent vectors is easy to compute this can give a useful alternative approach for the calculation of the geodesic distance between two unparametrized shapes $[q_0]$ and $[q_1]$: instead of minimizing in addition over the endpoint in the fiber of $[q_1]$ as we will describe in Section \ref{ssec:optimal_control_intrinsic}, one can alternatively minimize the horizontal energy between two fixed parametrized shapes. However, since, for the class of Sobolev metrics, calculating the horizontal projection is rather expensive, we will not follow this approach.

\subsection{The induced geodesic distance}\label{sec:geodesicdist}
Given a Riemannian metric one can consider the induced geodesic distance, which is defined as the infimum of the Riemannian length of all paths that connect two given points (shapes), i.e.,
\begin{equation}
d_{\Imm}(q_0,q_1) = \underset{q}{\operatorname{inf}}\; L(q):=\underset{q}{\operatorname{inf}} \int_0^1 \sqrt{G_{q(t)}(\partial_t q(t),\partial_t q(t))} dt\;,  
\end{equation}
where this infimum is calculated over all paths of immersions 
\begin{equation}
q:[0,1]\rightarrow \Imm_{C^{\infty}}(M,\R^d)
\end{equation}
with $q(0)=q_0$
and $q(1)=q_1$. We also obtain an induced geodesic distance function on the space of unparamtrized shapes:
\begin{equation}
d_{\mathcal S}([q_0],[q_1]) = \underset{\tau\in\Diff_{C^{\infty}}(M) }{\operatorname{inf}} 
d_{\Imm}(q_0,q_1\circ \tau)\;.
\end{equation}
In a similar way, one can obtain a geodesic distance function on spaces of shapes modulo translations and rotations. To factor out the group of scalings, one has to define a scale invariant metric on the space of parametrized shapes, cf. Remark~\ref{rem:invariances}.

On a finite-dimensional manifold the geodesic distance always defines a ``true'' metric, i.e.,
it is symmetric, satisfies the triangle inequality and non-degenerate. In infinite dimensions this result might not be true anymore as the geodesic distance can degenerate or even vanish identically. The first instance of this phenomenon has been found by Eliashberg and Polterovich for the $H^{-1}$-metric on the symplectomorphism group in \cite{eliashberg1993biinvariant}. In the context of shape analysis of submanifolds a similar result has been found by Michor and Mumford for the $L^2$-metric on the space of unparametrized submanifolds and on the diffeomorphism group, see \cite{Michor2005}. These results have been later extended to  spaces of parametrized submanifolds and to fractional-order metrics on diffeomorphism groups, see \cite{Bauer2013c,Bauer2012c}.  For Sobolev metrics of order $s\geq 1$ on spaces of submanifolds the following theorem shows that this ill-behavior cannot appear, which renders this class of metrics relevant for applications in shape analysis.   
\begin{theorem}[\cite{Michor2007,Bauer2011b}]
The geodesic distance of the Sobolev metric~\eqref{eq:RiemMetric} of order
$s\geq 1$ is non-degenerate both on the space of parametrized and unparametrized shapes.
\end{theorem}
This result has only been formally proven for the space of unparametrized shapes. For parametrized shapes it follows from the existence of simplifying transformation for $H^1$-metrics \cite{kurtek2012elastic,kurtek2010novel,Jermyn2011,Michor2008a} and the fact that $H^s$-metrics can be bounded from below by the $H^1$-metric for $s\geq 1$.

\subsection{The geodesic equation}
A classical result of Riemannian geometry states that the squared geodesic distance is also the infimum over all paths of the Riemannian energy,
\begin{equation} \label{eq:EnergyFunctional}
\begin{aligned}
d_{\Imm}(q_0,q_1)^2= \underset{q}{\operatorname{inf}} E(q) = \underset{q}{\operatorname{inf}} \int_0^1 G_{q(t)}(\partial_t q(t), \partial_t q(t)) dt \,.
\end{aligned}
\end{equation}
Thus we can find critical points of the length functional (i.e., locally shortest paths) by looking for critical points of the Riemannian energy functional. These critical points of the energy functional are called \emph{geodesics} and the first order condition for critical points, $dE(c) = 0$ is the \emph{geodesic equation}.
For our class of metrics the geodesic equation is a non-linear partial differential equation for the function $q = q(t,x)$, that is of order two in $t$ and of order $2s$ in $x\in M$. We refrain from presenting this equation here as it is rather complicated and not very insightful. Instead we refer the interested reader to the articles \cite{Bauer2011b,Michor2007}.

Since we are working in infinite dimensions the existence of geodesics is a highly non-trivial question and analytic solution formulas are only available in scarce special cases. 
The following theorem summarizes local and global well-posedness results, that are known for the class of reparametrization-invariant metrics on the space of parametrized submanifolds. 

\begin{theorem}
Let $M$ be a closed manifold and let $G$ be the reparametrization-invariant Sobolev metric of order $s \in \mathbb N$, $s\geq 1$, as defined in \eqref{eq:RiemMetric} on the space of parametrized submanifolds.   Then for any $r>\frac{d}{2}+1$ the initial  value  problem  for  the  geodesic  equation has  unique  local solutions  in  the  manifold $\Imm_{H^{r+2s}}(M,\mathbb R^d)$. The  solutions
depend  smoothly  on $t$
and  on  the  initial  conditions $q(0,\cdot)$
and $q_t(0,\cdot)$ and the
domain of existence (in $t$) is uniform in
$r$ and thus this also holds in the smooth category $\Imm_{C^{\infty}}(M,\mathbb R^d)$
\end{theorem}

Longtime existence of geodesics on the manifold $\Imm(M,\mathbb R^d)$ remains an open question for general parameter manifolds $M$. For the case $M=\S^1$, i.e., the space of curves in $\mathbb R^d$, the completeness of metrics of order $s\geq 2$ has been shown by Bruveris, Michor and Mumford \cite{Bruveris2014,Bruveris2014b_preprint}, see also \cite{nardi2016geodesics}.
\begin{theorem}\label{thm:completeness}
Let $G$ be the Sobolev metric of order $s \in \mathbb N$, $s\geq 2$, on the space of closed curves. We have:
\begin{enumerate}
    \item The Riemannian
manifold $\left(\Imm_{C^{\infty}}(\S^1,\R^d),G\right)$ is geodesically complete, i.e.,
given any initial conditions $(c_0,u_0)\in T \Imm_{C^{\infty}}(\S^1,\R^d)$ the solution of the geodesic 
equation for the metric $G$ with initial values $(c_0,u_0)$
exists for all time.
   \item  The metric completion of the space $\Imm_{C^{\infty}}(\S^1,\R^d)$ equipped with the geodesic distance 
   $\operatorname{dist}^{G}$ is the space $\Imm_{H^{s}}(\S^1,\R^d)$ of immersions of Sobolev class $H^s$.
   Furthermore, any two curves in the same connected component of the space $\Imm_{H^{s}}(\S^1,\R^d)$ can be joined by a minimizing geodesic.
\end{enumerate}
\end{theorem}
There is some evidence that suggests that a similar result should be valid for general $M$ if the order $s$ of the metric satisfies $s>\frac{d}2+1$. A  proof of this statement is, however, largely missing.

\subsection{An optimal control formulation of the geodesic problem on the space of unparametrized shapes}
\label{ssec:optimal_control_intrinsic}
In the following we will describe the optimal control formulation for the geodesic boundary value problem, as introduced in Section~\ref{sec:control}, in the specific situation studied in this Section. We recall, that 
a geodesic between $b_0 = \pi(q_0)=[q_0]$ and $b_1 =\pi(q_1) = [q_1]$ minimizes $\int_0^1 G_{q(t)}(h(t),h(t)) dt$ subject to $\dot{q}(t) = h(t) \in T_{q(t)} \Imm(M,\R^d)$ and the boundary constraints $q(0) = q_0$ and $q(1) \in \pi^{-1}(b_1)$. The latter constraint can  be replaced by jointly optimizing over variations in the fiber $q(1)\circ \tau$ for $\psi\in \Diff(M)$. This leads to the following formulation of the boundary value problem on the space of unparametrized shapes:
\begin{equation}
\label{eq:exactmathcing_intrinsic}
    \min_{h,\tau} \int_0^1 G_{q(t)}(h(t),h(t)) dt
\end{equation}
over $h \in L^2([0,1],C^{\infty}(M,\R^d))$, $\tau \in \Diff(M)$ with $q(0)=q_0$, $q(1)=q_1 \circ \tau$ and $\dot{q}(t) = h(t)$. It is often more practical to relax the end-time constraint using a fidelity term blind to fiber variations, in other words to reparametrizations of $q(1)$ and $q_1$. This is where the chordal distances as defined in Section \ref{sec:chordal} fit in. In particular, using oriented varifold distances, the inexact matching problem may be formulated as:   
\begin{equation}
\label{eq:inexactmathcing_intrinsic}
    \min_{h} \int_0^1 G_{q(t)}(h(t),h(t)) \ dt + \lambda d_{\text{var}}(q(1),q_1)^2
\end{equation}
with $\dot{q}(t) = h(t)$, $q(0)=q_0$ and $\lambda>0$ being a weight parameter. Interestingly, \eqref{eq:inexactmathcing_intrinsic} turns into an infinite-dimensional \textit{optimal control} problem in which, at each time $t$, the state of the system is the immersion $q(t)$ and the control is the vector field $h(t)$.

\subsection{Numerical aspects}
We will now describe a discretization of the relaxed geodesic estimation problem \eqref{eq:inexactmathcing_intrinsic}. We will focus on the discretization of the intrinsic metric part of the cost, the  discretization of the varifold metric having been discussed in Section~\ref{sec:chordal}. In the language of Section~\ref{sec:control}, our approach will be based on a trajectory optimization approach. For the space of curves and Sobolev metrics of first order, there exist simplifying transformations --related to the so-called square root velocity mapping-- that map the manifold of shapes into a subset of a flat space. This in turn allows one to obtain (almost) explicit formulas for geodesics and geodesic distance, which makes the numerical computations on the space of parametrized shapes trivial \cite{Younes1998,Jermyn2011,Michor2008a,Bauer2014b,KuNe2018}. Since we aim to describe a numerical framework for a more general class of metrics, we will not describe this approach, but focus on the discretization of the energy functional.  Furthermore, we will focus on the case of open and closed curves. For a four-parameter family of second order Sobolev metrics the presented framework is available under an open source license\footnote{The code can be downloaded at github: 
\url{https://www.github.com/h2metrics/h2metrics}, see \cite{BauerCoCV2018} for a detailed description of the implementation.}.

Theoretically the described approach will directly generalize to higher dimensional objects. However, so far only the case of open and closed curves has been implemented for this general class of metrics and is available in open access.  It is an ongoing project of the authors to obtain a similar numerical framework for the space of surfaces in $\mathbb R^3$. (Note, however, that, based on a generalization of the square root velocity transform, a numerical scheme for a particular Sobolev metric of first order on the space of surfaces has been developed in \cite{laga2017numerical,kurtek2013landmark}.)\\

To discretize the energy of a path of curves with respect to the intrinsic metric $G$ it is convenient to use a B-splines representation of the involved elements, i.e., a path of curves $c(t,\theta)$ is represented by a tensor-product B-spline on knot sequences of orders $n_t$ in time and $n_\theta$ in space. Note that the degree of the spline in space $n_\theta$ should be greater than the order of the metric $s$. Using this discretization one obtains $N_t N_\theta$ basis splines $B_i(t)C_j(\theta)$, where $i\in \{1,\ldots,N_t\}$ and $j\in \{1,\ldots,N_\theta\}$ and $N_t$ and $N_\theta$ are the number of control points in each variable respectively. Note that $B_i(t)$ are B-splines defined by an equidistant simple knot sequence on $[0,1]$ with full multiplicity at the boundary knots, and $C_j(\theta)$ are defined by an equidistant simple knot sequence on $[0,2\pi]$ with either periodic boundary conditions or full multiplicity at the boundary for closed or open curves respectively.
Now, paths of curves $c(t,\theta)$ are parametrized as
\begin{equation}
\label{eq:TensorPath}
c(t, \th) = \sum_{i=1}^{N_t} \sum_{j=1}^{N_\th} c_{i,j} B_i(t) C_j(\th)\,.
\end{equation}
More details on this approach can be found  in~\cite{BauerCoCV2018,BBHM2017}. The full multiplicity of the boundary knots in $t$ simplify the boundary constraints, since we have
\begin{align*}
c(0, \th) &= \sum_{j=1}^{N_\th} c_{1,j} C_j(\th)\,, &
c(1, \th) &= \sum_{j=1}^{N_\th} c_{N_t,j} C_j(\th)\,.
\end{align*}
Thus the initial curve $c(0)$ is given by the control points $c_{1,j}$ only. 
Using this B-spline representation, it is then straightforward to discretize the Riemannian energy: we specifically use Gaussian quadrature with quadrature nodes placed between knots where the curves are smooth. 

This discretization then yields a trajectory optimization approach on the spline control points, which in that case boils down to an unconstrained minimization problem --the free variables being the intermediate control points $c_{i,j}$ for $i=2,\ldots,N_{t-1}$-- which can be tackled with standard methods of finite-dimensional optimization such as the L-BFGS algorithm from the HANSO library \cite{hanso} that we use. The complexity of the algorithm at each optimization step is approximately $O(N_{\theta}N_t)$ (evaluation of splines and their derivatives to compute the  energy and its gradient for which the required number of basic operations is linear in the total number of spline control points). To solve the geodesic matching problem, one also needs to compute at each optimization step the varifold distance and its gradient which have complexity $O(N_{\theta}^2)$.

\begin{figure}
    \centering
 
        \includegraphics[trim={4cm 0 3.5cm 0},clip,width=\textwidth]{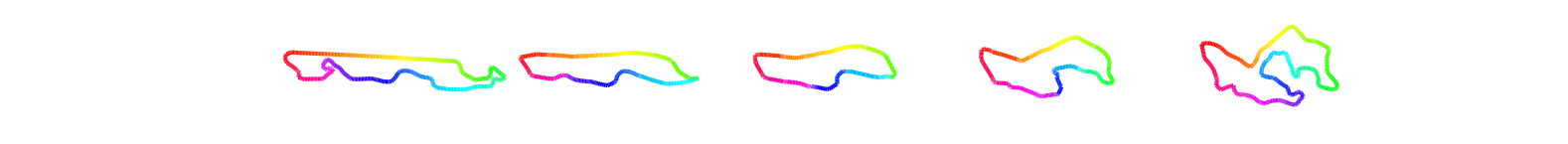}
     \includegraphics[trim={4cm 0 3.5cm 0},clip,width=\textwidth]{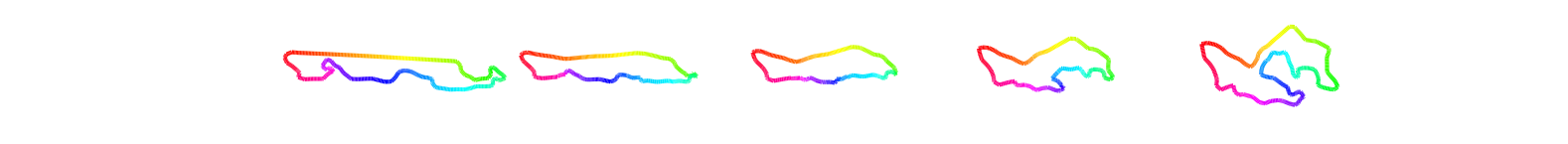}

    \caption{Optimal deformations between two unparametrized shapes with respect to an intrinsic $H^1$-metric at selected time steps ($t=0, 0.25, 0.5, 0.75, 1$). These boundary shapes represent the outlines of Florida 5th congressional district and Maryland 7th congressional district. In the second row the amount of $H^1$-energy is increased by a factor of 100. One can observe that the higher curvature regions are easier resolved in the first row as compared to the second row.}
    \label{fig:geodesic:intrinsic}
\end{figure}

We show examples of geodesics on the space of unparametrized curves in Figure~\ref{fig:geodesic:intrinsic}. On can clearly see the influence of the choice of metric constants on the resulting geodesics.  As a comparison we will present the same example for the outer and hybrid deformation metric model in Figure~\ref{fig:outer_hybrid}.

\section{Outer deformation metric models}
\label{sec:LDDMM}

We now switch to a construction based on the transitive group action model described in Section \ref{sec:submersion} paragraph (2), in which $\Imm(M, \R^d)$ becomes the base space and the top space is now provided by the diffeomorphism group of  $\R^d$ ($\Diff(\R^d)$) with action $\phi\cdot q = \phi\circ q$.

If $\bB$ is a Banach space of vector fields over $\R^d$, we let $\Diff_\bB$ denote the space of diffeomorphisms $\phi$ on $\R^d$ such that $\phi-\id$ and $\phi^{-1}-\id$ both belong to $\bB$. We will in particular use $\Diff_{C_0^r}(\R^d)$ to denote the space corresponding to 
 $\bB=C_0^r(\R^d, \R^d)$, the space of vector fields $v$ on $\R^d$ that are $r$ times continuously differentiable (with $r\geq 1$) and tend to 0 (together with their first $r$ derivatives) uniformly at infinity. Such spaces are Banach manifolds with $\bB = T_\phi \Diff_\bB$ for all $\phi$ (and in particular at $\phi=\id$). We will always assume that $\bB$ is continuously embedded in $C_0^1(\R^d, \R^d)$.
 
 In this infinite-dimensional setting, the Riemannian framework needs to be weakened into a sub-Riemannian one, in which the metric will only be defined on subspaces of the tangent spaces, inducing the notion of admissible trajectories. More precisely, let $\bV$ be a Hilbert  
of vector fields, continuously embedded in $\bB$, and let $\scp{\cdot}{\cdot}_\bV$ and $\|\cdot\|_\bV$ denote the inner product and norm. For $\phi\in\Diff_\bB$, we define
\[
\bV_\phi = \bV\circ \phi = \{v\circ \phi, v\in \bV\}
\]
(so that $\bV = \bV_\id$), and the sub-Riemannian metric
\[
G_\phi(w, w') = \scp{w\circ \phi^{-1}}{w'\circ\phi^{-1}}_\bV, \text{ for } w, w'\in \bV_\phi. 
\]
(This metric is called sub-Riemannian because it is only defined on a subset of the tangent bundle of $\Diff_\bB$.)

A differentiable path $\psi(\cdot)$  on $\Diff_\bB$ is then called admissible if $\dot\psi(t) \in V_{\psi(t)}$ for all $t$ (such paths are often called horizontal in the sub-Riemannian geometry literature, but we are using horizontal with another meaning here). The energy of an admissible path (defined for $t\in [0,1]$) is given by
\[
E(\psi) = \int_0^1 G_{\psi(t)}(\dot\psi(t), \dot\psi(t))\, dt
\]
A diffeomorphism $\phi\in \Diff_\bB$ is called attainable (from the identity) is there exists an admissible path with finite energy such that $\psi(0) = \id$ and $\psi(1) = \phi$. The set of attainable diffeomorphisms will be called $\Diff_\bV^{\att}$. One also defines a right-invariant sub-Riemannian distance $d_V$ on $\Diff_\bV^\att$ such that $d_\bV(\phi, \phi') = d_\bV(\id, \phi'\circ\phi^{-1})$ is the square root of the minimal energy among all admissible paths between $\id$ and $\phi'\circ \phi^{-1}$. 

We then have the following result (see \cite{younes2010shapes} for a proof).
\begin{theorem}[Trouv\'e]
If $\bB$ is embedded in $C_0^r(\R^d, \R^d)$, then
$\Diff^\att_\bV$ is a subgroup of $\Diff_{C_0^r}(\R^d)$ and a complete space for the sub-Riemannian distance $d_\bV$. Moreover, for any $\phi\in \Diff_\bV^\att$, there exists an optimal path $\psi(\cdot)$ such that $E(\psi)$ is minimal over all admissible paths between $\id$ and $\phi$. 
\end{theorem}
Notice that $\Diff_\bV^\att$ is in general distinct from $\Diff_\bV$. Exceptions in which the two sets coincide are provided when $\bV$ equals a Sobolev space of vector fields of high enough order (see \cite{bruveris2014completeness}).  \\

$\Diff_\bB$ acts on both $\Imm_{C^s}(M, \R^d)$ and $\Emb_{C^s}(M, \R^d)$ as soon as $r\geq s$ through the action $\phi\cdot q = \phi\circ q$ (because, if $q$ is in any of these two spaces, $\phi\cdot q$ also belongs to the same space).  The infinitesimal action is $v\cdot q = v\circ q$. The optimal control problems described in equations \eqref{eq:r.sub.geod} and \eqref{eq:r.sub.geod.2} can then be rewritten as follows. The square distance on the base space between $q_0$ and $q_1$ is obtained by minimizing
\begin{equation}
\label{eq:lddmm.exact}
\int_0^1 \|v(t)\|^2_\bV\, dt
\end{equation}
subject to $q(0) = q_0$, $q(1) = q_1$ and $\dot q = v\circ q$. The relaxed version removes the constraint $q(1) = q_1$ and replaces it by a penalty $U(q(1), q_1)$ added to the integral.

Horizontal and vertical subspaces can be defined within this sub-Riemannian framework as subsets of $\bV$ or $\bV_\phi$, letting, for a given $q$, $\mf v_q = \{v \in \bV: v\circ q = 0\}$, $\mf h_q = \mf v_q^\perp$
with $V_\phi = \mf v_{\phi\circ q}\circ \phi$, $H_\phi = \mf h_{\phi\circ q}\circ \phi$. The specific instances taken by these spaces can however be surprising in some cases. It is for instance possible to find examples for which $\mf v_q=\{0\}$ so that $\mf h_q = \bV$, which is counter-intuitive compared to the finite dimensional setting in which $\B$ and $\mf h_q$ should have the same dimension.\\

As described in Section \ref{sec:submersion}, the induced (sub-Riemannian) metric on $\Imm_{C^s}(M, \R^d)$ is defined by
\begin{equation}
\label{eq:induced_metric_LDMMM}
G^\bV_q(h,h) = \min\{\|v\|_\bV^2: v\in \bV, v\circ q = h \}
\end{equation}
for all $h\in \bV_q = \bV\circ q\subset C^s(M, \R^d)$. If $h\in \bV_q$, then there is a unique $v = h^q\in \mf h^q$ such that $G^\bV_q(h,h) = \|h^q\|_\bV$, and the minimization in \eqref{eq:lddmm.exact} (or its relaxed version) can be reduced with the additional constraint $v(t) \in \mf h_{q(t)}$ without loss of generality. This remark will be useful, in particular, for the discrete version of the problem that is described below.\\

In order to allow this metric to formally descend to unparametrized curves $[q]\in \mathcal S$, we need to check that it is right invariant by diffeomorphisms $\tau\in \Diff_{C^s}(M)$, i.e., that 
\[
G^\bV_{q\circ \tau}(h\circ\tau ,h\circ \tau) = G^\bV_q(h,h)
\]
with $\bV_{q\circ \tau} = \bV_q\circ \tau$, which is clear from the definitions. So, if $U$ is any function such that $U(q, q') = 0$ if and only if $q'$ is a reparametrization of $q$ (which holds for the chordal metrics introduced in Section \ref{sec:chordal} for instance), the minimization of
\[
\int_0^1 \|v(t)\|^2_\bV\,dt + U(q(1), q_1)
\]
subject to $q(0) = q_0$ and $\dot q(t) = v(t)\circ q(t)$ provides a metric registration method between unparametrized immersions, inducing a method called large deformation diffeomorphic metric mapping, or ``LDDMM'' \cite{Beg2005,Glaunes2006,joshi2000landmark,younes2010shapes,Glaunes2008}.  This problem can also be expressed in  sub-Riemannian form, minimizing
\[
\int_0^1 G^\bV_{q(t)}(\dot q(t), \dot q(t))\,dt + U(q(1), q_1)
\]
subject to $q(0) = 0$, $\dot q \in \bV_{q(t)}$ for all $t$. Notice that any admissible path starting from $\Emb_{C^s}(M, \R^d)$ will remain in that space at all times, and that, when restricted to that space, this optimal control problem becomes one between embedded submanifold of $\R^d$ (diffeomorphic to $M$).
\\

Unsurprisingly, the complications that result from the infinite-dimensional nature of the shape space disappear after discretization, in which case the variational framework is a direct application of the formulation described in Section \ref{sec:submersion}. Using discrete curves and triangulated surfaces, the manifold $M$ can be replaced by a simplicial complex, which, in the surface case, is provided by a family of singletons, pairs and triples of integers. In the special case of a triangulated surface, the complex can be described as a finite set of indices, say $\I$ together with a set $\T$ consisting of triples $(i,j,k)$ where triples that differ with a circular permutation are identified, i.e., $(i,j,k) = (j,k,i)$. Edges are deduced from $\I$ and $\T$ as the set containing all pairs that are subsets of triples in $\T$. Letting $\M = (\I, \T)$, We can then denote by $\Imm(\M, \R^3)$ the set of all one-to-one mappings $q = (q(i), i\in \I)$ from $\I$ to $\R^3$, and $(q, \T)$ then forms a triangulated surface. (Typically, $\M$ is obtained though a triangulation of a continuous surface $M$.) With the action of diffeomorphisms still defined by $\phi\cdot q = \phi\circ q$, and the infinitesimal action $v\cdot q = v\circ q$, the vertical space $\mf v_q$ is given by
\[
\mf v_q = \{v\in \bV: v(q(i)) = 0, i\in\I\}.
\]The space $\mf h_q = \mf v_q^\perp$ can be described using the reproducing kernel of $\bV$, which is an RKHS of vector fields since it is embedded in $C^1_0(\R^r, \R^d)$. Denoting this kernel by $K$, which is therefore a matrix-valued function $(x,y) \mapsto K(x,y)$, the space $\mf h_q$ is simply
\[
\mf h_q = \left\{v = \sum_{i\in\I} K(\cdot, q(i)) a(i): a(i)\in\R^3, i\in\I\right\}.
\]
Moreover, if $v = \sum_{i\in\I} K(\cdot, q(i)) a(i)\in \mf h_q$, then 
\[
\|v\|_\bV^2 = \sum_{i,j\in\I} a(i)^TK(q(i),q(i))a(i)
\]
is also explicit. This implies that the problem in \eqref{eq:lddmm.exact} with the additional constraint $v\in \mf h_q$ can now be rephrased as a finite dimensional optimal control problem, minimizing
\[
\int_0^1 L(q(t), a(t))\, dt
\]
subject to $q(0, \cdot) = q_0$, $q(1, \cdot) = q_1$ and $\dot q(t,i) = \sum_{j\in\I} K(q(t,i), q(t,j)) a(t,j)$, where
\[
L(q,a) = \sum_{i,j\in\I} a(i)^TK(q(i),q(j))a(j). 
\]
The corresponding metric on $\Imm(\M, \R^3)$ is actually Riemannian (not just sub-Riemannian), with 
\[
G^\bV_q(h,h) = h^T \K(q)^{-1} h
\]
where $h = (h(i), i\in\I)$ with $h(i) \in \R^3$ is identified with a column vector and $\K(q)$ is the matrix formed with blocks $K(q(i), q(j))$, $i,j\in\I$.

Finally, the relaxed version of this problem requires using data attachment terms $U(q, q')$ adapted to triangulated surfaces, such as those described in Section \ref{sec:chordal}. Several algorithms and open access numerical implementations of curve/surface LDDMM have been developed. In the vocabulary of Section \ref{sec:control}, shooting-based implementations include codes such as Deformetrica\footnote{\url{http://www.deformetrica.org/}}, fshapesTk\footnote{\url{https://github.com/fshapes/fshapesTk}}, catchine\footnote{\url{www.mi.parisdescartes.fr/~glaunes/catchine.zip}} while other implementations rely instead on trajectory optimization over the time-varying costates $a(t)$, such as py-lddmm\footnote{\url{https://bitbucket.org/laurent_younes/registration.git}}. Matching results based on the latter code are shown in the first rows of Figure \ref{fig:outer_hybrid} for curves and of Figure \ref{fig:outer_hybrid_surfaces} for surfaces.

\section{A hybrid metric model}
\label{sec:hybrid}
The outer metric framework of the previous section has the advantage of enforcing a diffeomorphic transformation of the template shape and prevents geodesics from changing topology along optimal paths, which is often an important constraint in  applications. Yet, penalizing an extrinsic deformation of the ambient space as measured by a global metric $G_q^{V}$  can also end up being limited. This is the case, most notably, when each shape is itself a complex of multiple curves or surfaces in close proximity to one another \cite{Durrleman2014,Arguillere2016}. Indeed, while one expects smooth deformations of each individual component, it is usually critical to also estimate more local changes of alignment and sliding motions between them: these two conditions are generally hard to jointly satisfy in the LDDMM model, even with multiscale approaches. In this section, we take the advantages of both intrinsic and outer metrics and combine them into an hybrid model. The following is a condensed version of the more complete presentation of \cite{Younes2018_hybrid}. 

The construction of hybrid metrics follows a similar setting as Section \ref{sec:LDDMM}. The key difference consists in replacing the cost function in \eqref{eq:lddmm.exact} by:
\begin{equation}
\label{eq:hybrid.exact}
\int_0^1 \|v(t)\|^2_{q(t)}\, dt
\end{equation}
where, for $q \in \Imm_{C^s}(M, \R^d)$,  $\|v\|_{q}^2 = \|v\|_{\bV}^2 + G_{q}(v\circ q,v\circ q)$, $G_q$ being an intrinsic metric of order at most $s$ such as the Sobolev metric of \eqref{eq:RiemMetric}. Note that the LDDMM framework corresponds to the particular case $G_q =0$. Now the induced metric on $\Imm_{C^s}(M, \R^d)$ in \eqref{eq:induced_metric_LDMMM} becomes
\begin{equation}
\label{eq:induced_metric_hybrid}
G^\bV_q(h,h) = \min\{\|v\|_\bV^2 + G_{q}(h,h): v\in \bV, v\circ q = h \}
\end{equation}
for all $h\in \bV_q = \bV\circ q\subset C^s(M, \R^d)$. Note that the two terms in this hybrid model play a different and complementary role. The intrinsic term constrains the regularity of the deformation of the shape $q(M)$ itself while the LDDMM term is mostly here to enforce the global transformation resulting from the full vector field $v$ to be diffeomorphic. In particular, small scale kernels for the norm on $\bV$ can be typically used in this context, in order, for example, to allow more local transformations of the background space.  

The inexact matching problem under the hybrid metrics then consists in the minimization of:
\begin{equation*}
    \int_0^1 \left[ \|v(t)\|_{\bV}^2 + G_{q(t)}(v(t)\circ q(t),v(t)\circ q(t)) \right] \,dt + F(q(1), q_1)
\end{equation*}
over $v \in L^2([0,1],\bV)$ and subject to $q(0) = q_0$, $\dot{q}(t) = v(t) \circ q(t)$. As in Section \ref{sec:LDDMM}, if one considers discrete simplicial shapes, this problem can be turned into a finite-dimensional optimal control problem. Specifically, using the same notations as previously, the discrete Lagrangian is now:
\begin{equation*}
    L(q,a) =  a^T \K(q)a + (\K(q) a)^T \Lambda(q) (\K(q) a)
\end{equation*}
where $\K(q)$ is again the matrix of the kernel of $V$ evaluated at the vertices, $\K(q) a \in \R^{d |\I|}$ the vectorized deformation field at all vertices, and  $\Lambda(q) \in \R^{(d|\I|) \times (d|\I|)}$ is the matrix that corresponds to the discretization of the differential operator defining the intrinsic part of the metric. $\Lambda(q)$ is typically a sparse matrix obtained from either a finite difference or a finite element scheme. For instance, \cite{Younes2018_hybrid} uses an $H^1$ metric for $G_q$ which is discretized based on a finite difference scheme on each simplex of the shape. In that case, the resulting approximations $(\K(q) a)^T \Lambda(q) (\K(q) a)$ become:
\begin{equation}
\label{eq:curve.h1}
    \sum_{(i,j) \in \T} \frac{|h(j)-h(i)|^2}{|q(j)-q(i)|}  \ \ \ \text{(curve)}
\end{equation}
\begin{equation}
\label{eq:surf.h1}
    \sum_{(i,j,k) \in \T} \frac{1}{4A_{ijk}} |h(i) (q(k)-q(j)) + h(j) (q(i)-q(k)) + h(k) (q(j)-q(i))|^2 \ \ \ \text{(surface)}
\end{equation}
where $h \in \R^{|\I| \times d}$ are the values of $v$ at the different vertices $q(i)$ and, in the case of triangulated surfaces, for each triangle $(i,j,k)$ in $\T$, $A_{ijk}$ denotes its area which is given by the norm of $(q(j)-q(i))\wedge(q(k)-q(i))$.

\begin{figure}
\includegraphics[trim=5cm 5cm 4.5cm 5cm, clip,width=0.19\textwidth]{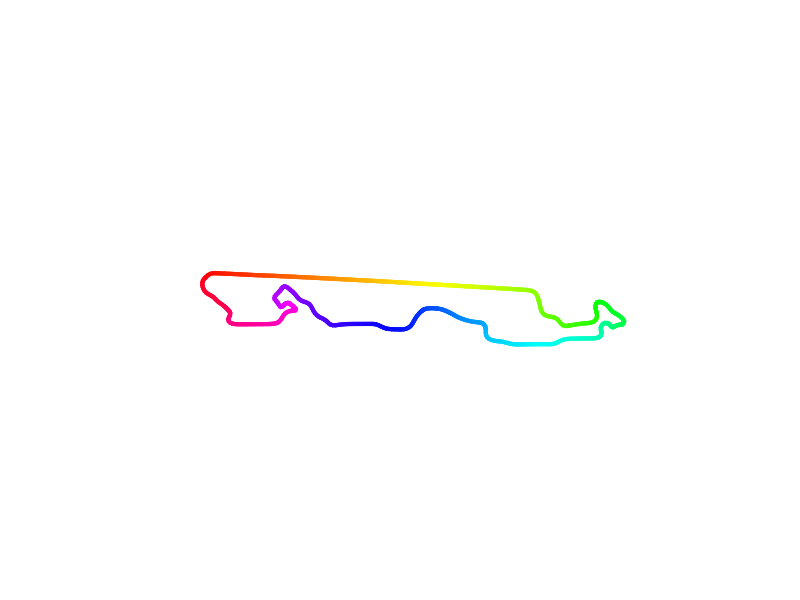}
\includegraphics[trim=5cm 5cm 4.5cm 5cm, clip,width=0.19\textwidth]{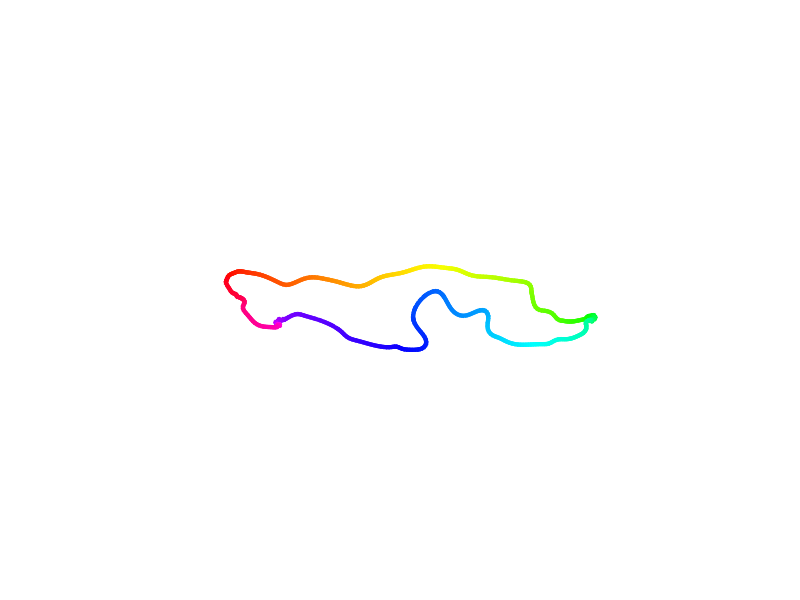}
\includegraphics[trim=5cm 5cm 4.5cm 5cm, clip,width=0.19\textwidth]{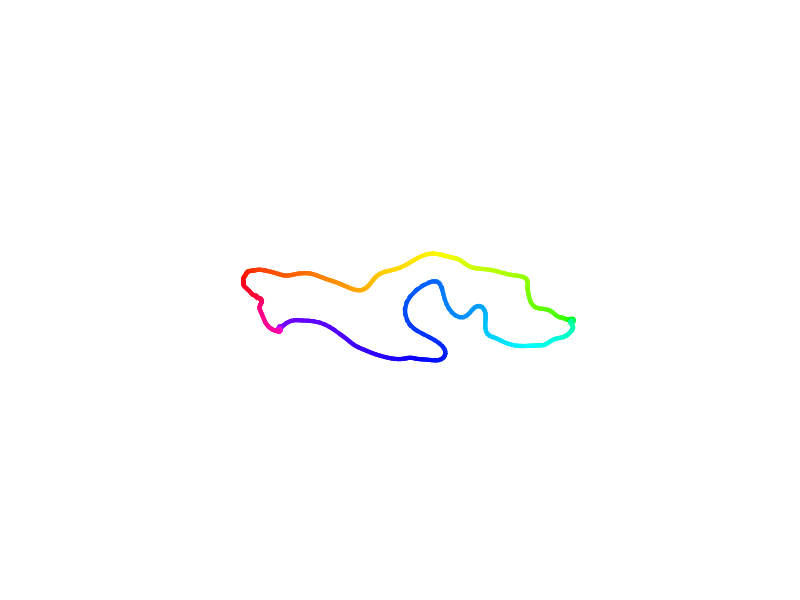}
\includegraphics[trim=5cm 5cm 4.5cm 5cm, clip,width=0.19\textwidth]{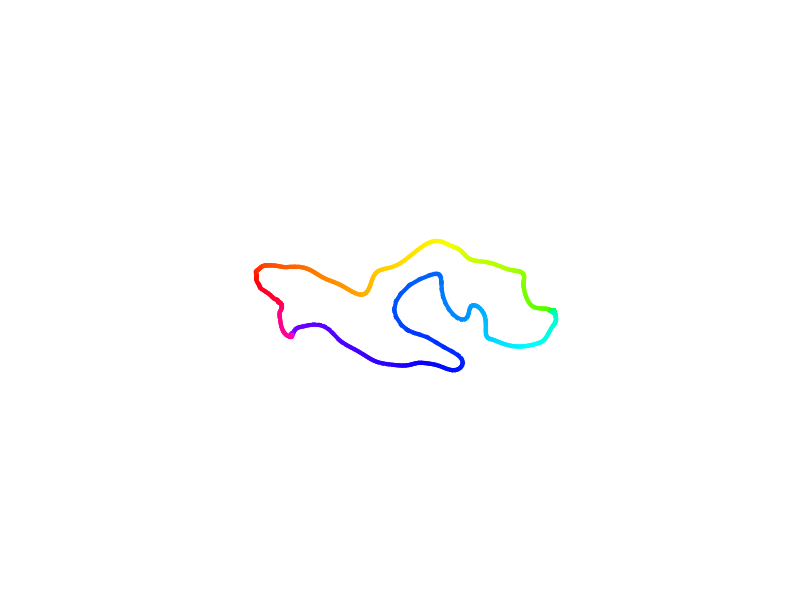}
\includegraphics[trim=5cm 5cm 4.5cm 5cm, clip,width=0.19\textwidth]{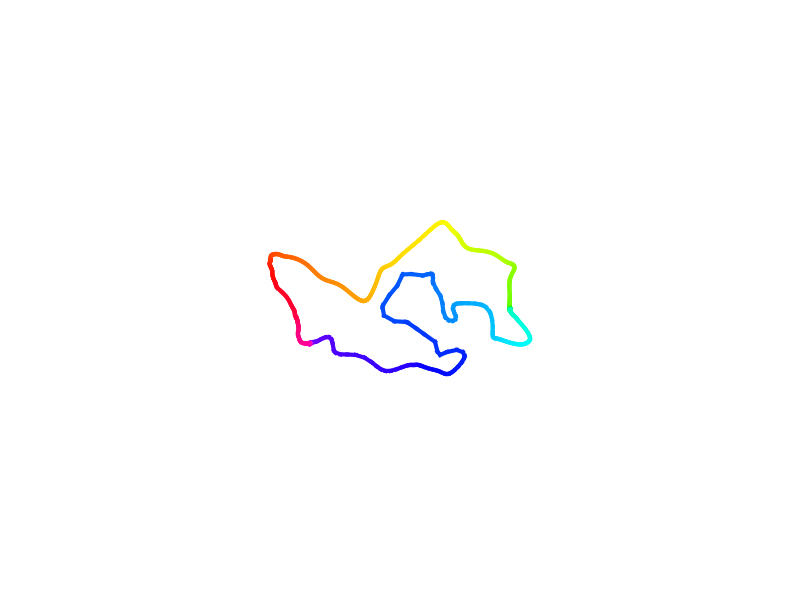}\\

\includegraphics[trim=5cm 5cm 4.5cm 5cm, clip,width=0.19\textwidth]{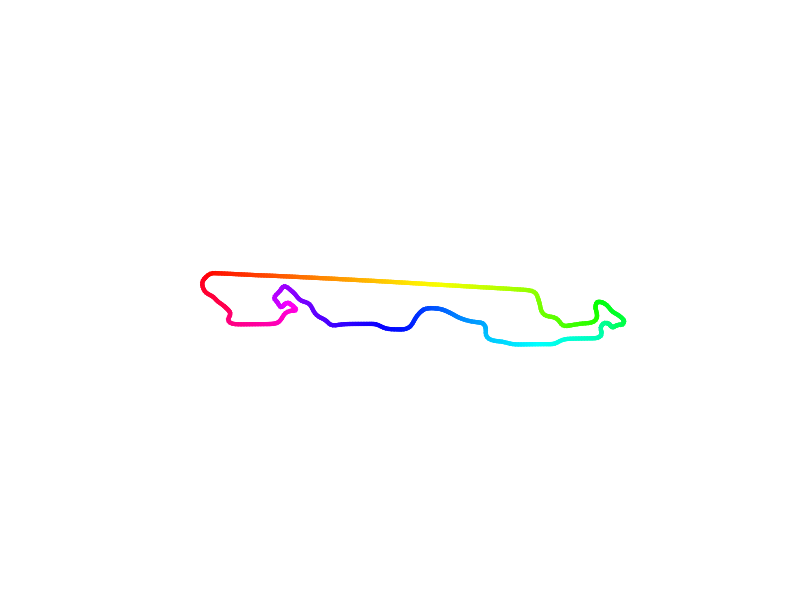}
\includegraphics[trim=5cm 5cm 4.5cm 5cm, clip,width=0.19\textwidth]{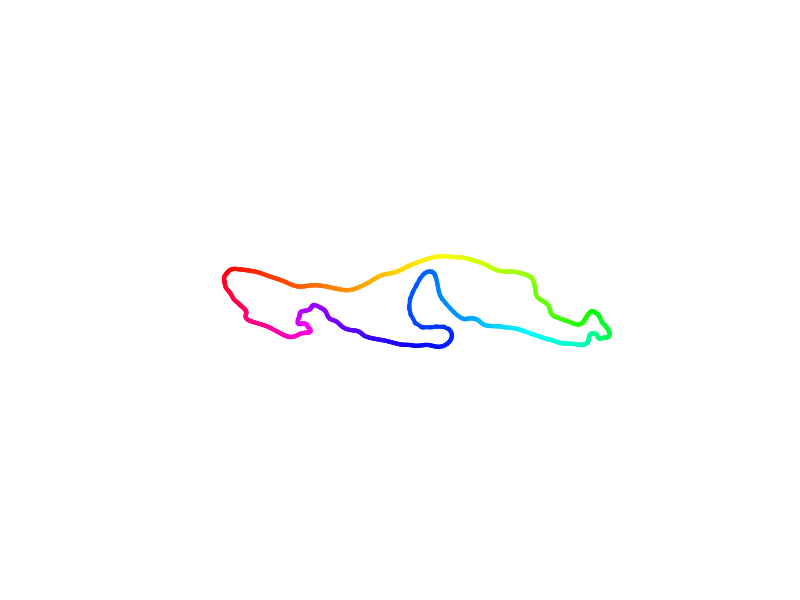}
\includegraphics[trim=5cm 5cm 4.5cm 5cm, clip,width=0.19\textwidth]{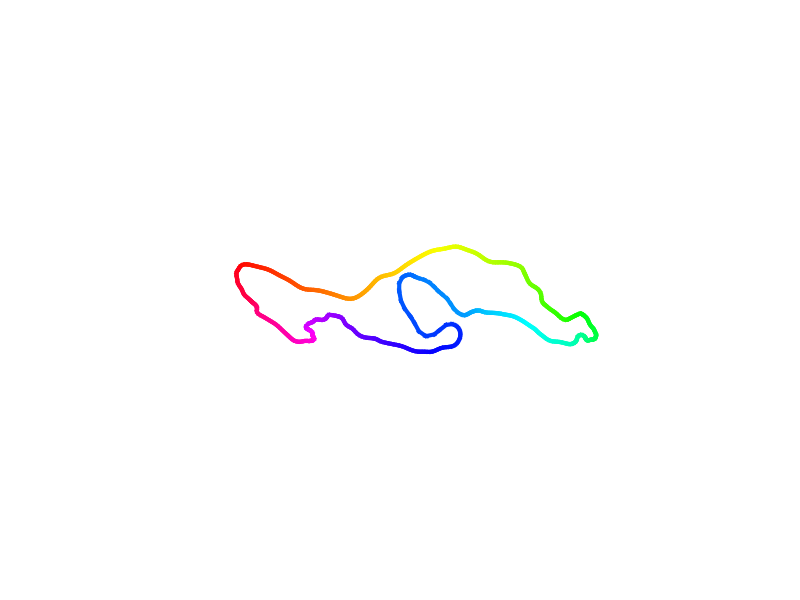}
\includegraphics[trim=5cm 5cm 4.5cm 5cm, clip,width=0.19\textwidth]{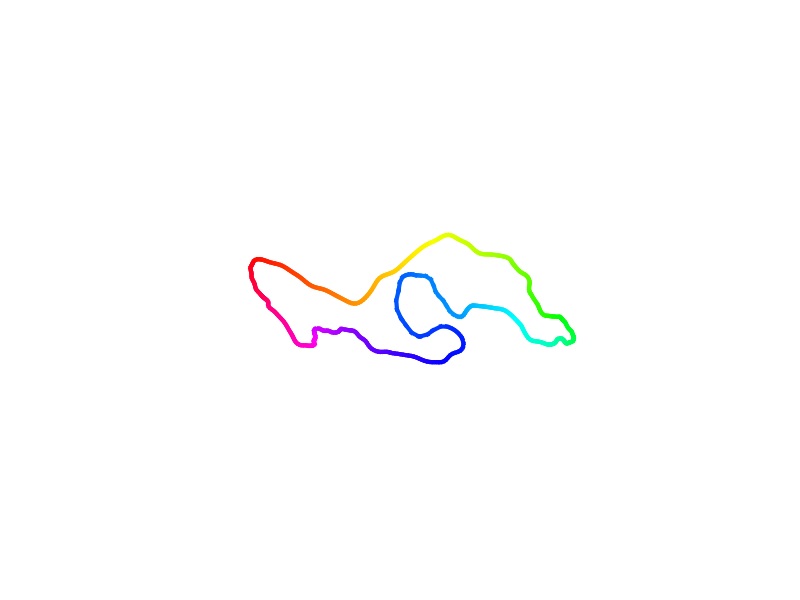}
\includegraphics[trim=5cm 5cm 4.5cm 5cm, clip,width=0.19\textwidth]{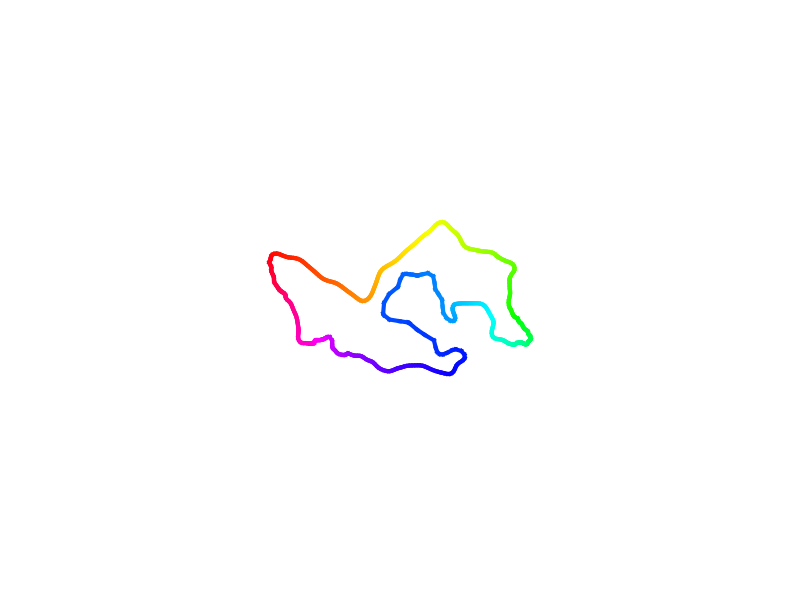}\\

\caption{
\label{fig:outer_hybrid} Geodesics for outer and hybrid metrics on closed curves at selected time steps ($t=0, 0.3, 0.5, 0.7, 1$). Top row (Outer metric): a few portions of the template curve are dramatically compressed to align the two shapes (in the green region on the right and in the purple/pink transition), at a point making them almost invisible at the image resolution. Bottom row (hybrid metric): the regions that were compressed in the outer-metric geodesic now properly unwrap since compression along the curve is now more strongly penalized. }
\end{figure}

\begin{figure}
    \centering
    \includegraphics[width=0.35\textwidth]{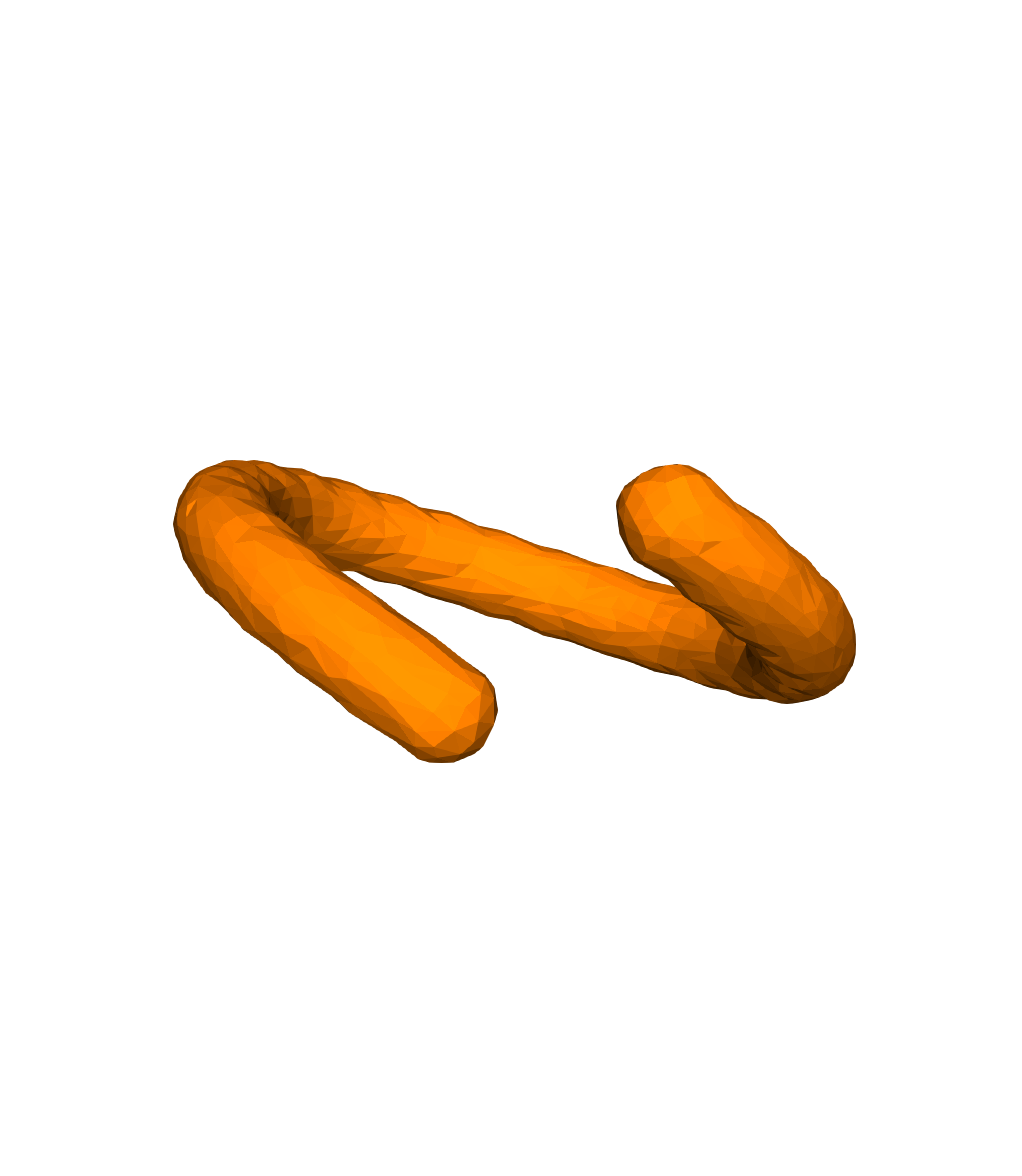}
    \includegraphics[width=0.35\textwidth]{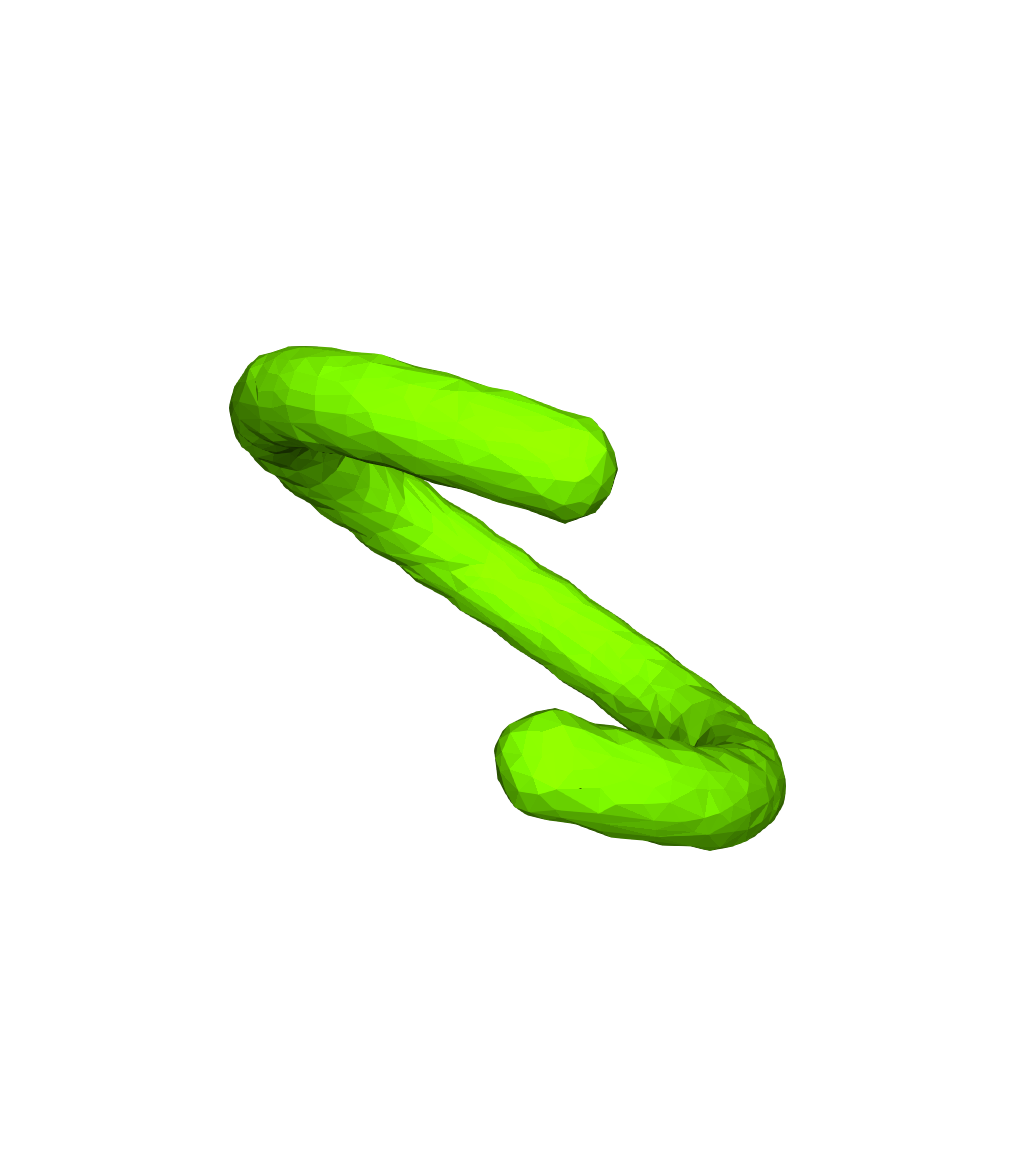}
    \caption{Deformed incomplete rings. Left: template; Right: target.}
    \label{fig:surf}
\end{figure}

\begin{figure}
\includegraphics[trim=5cm 5cm 4.5cm 5cm, clip,width=0.19\textwidth]{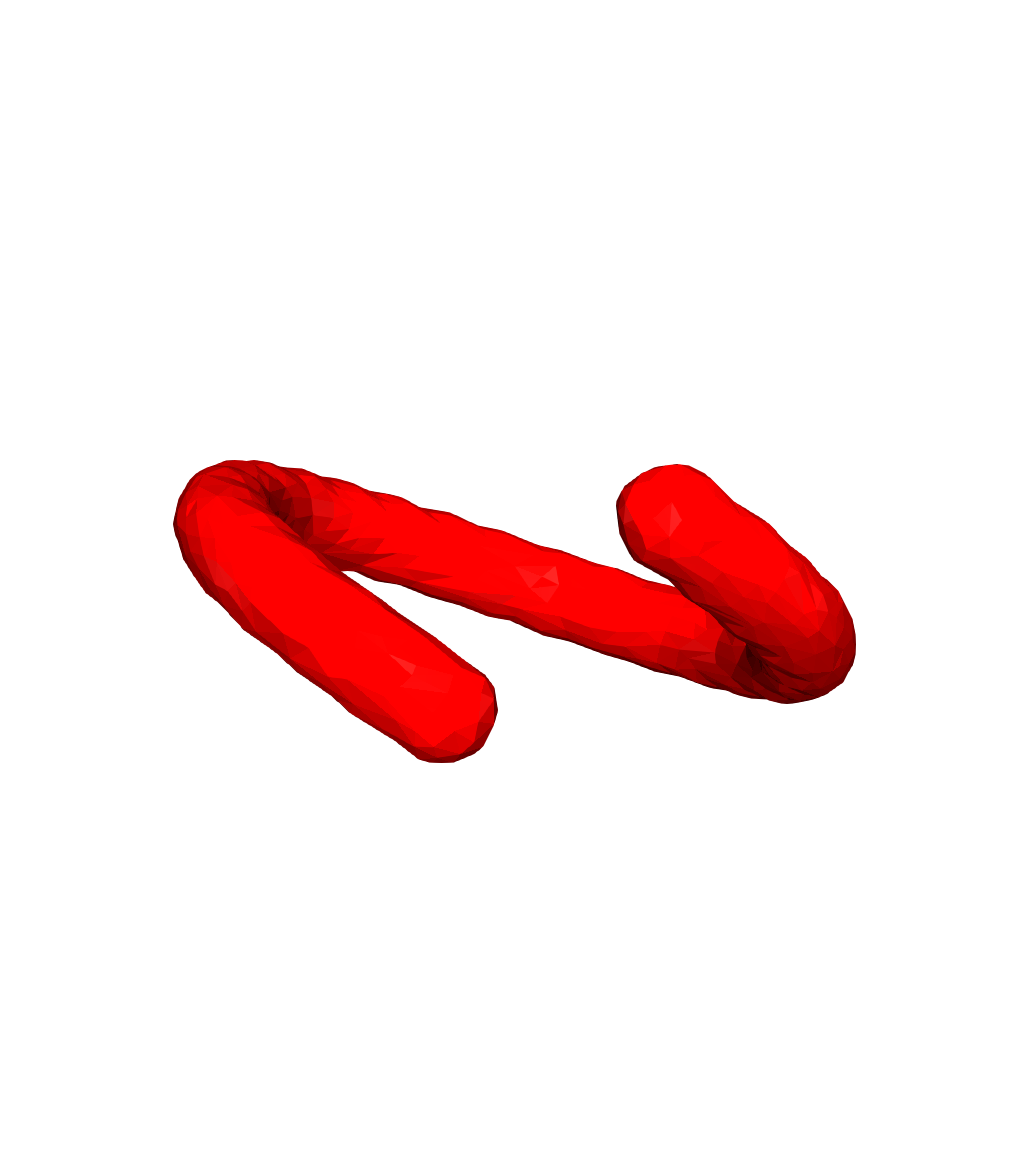}
\includegraphics[trim=5cm 5cm 4.5cm 5cm, clip,width=0.19\textwidth]{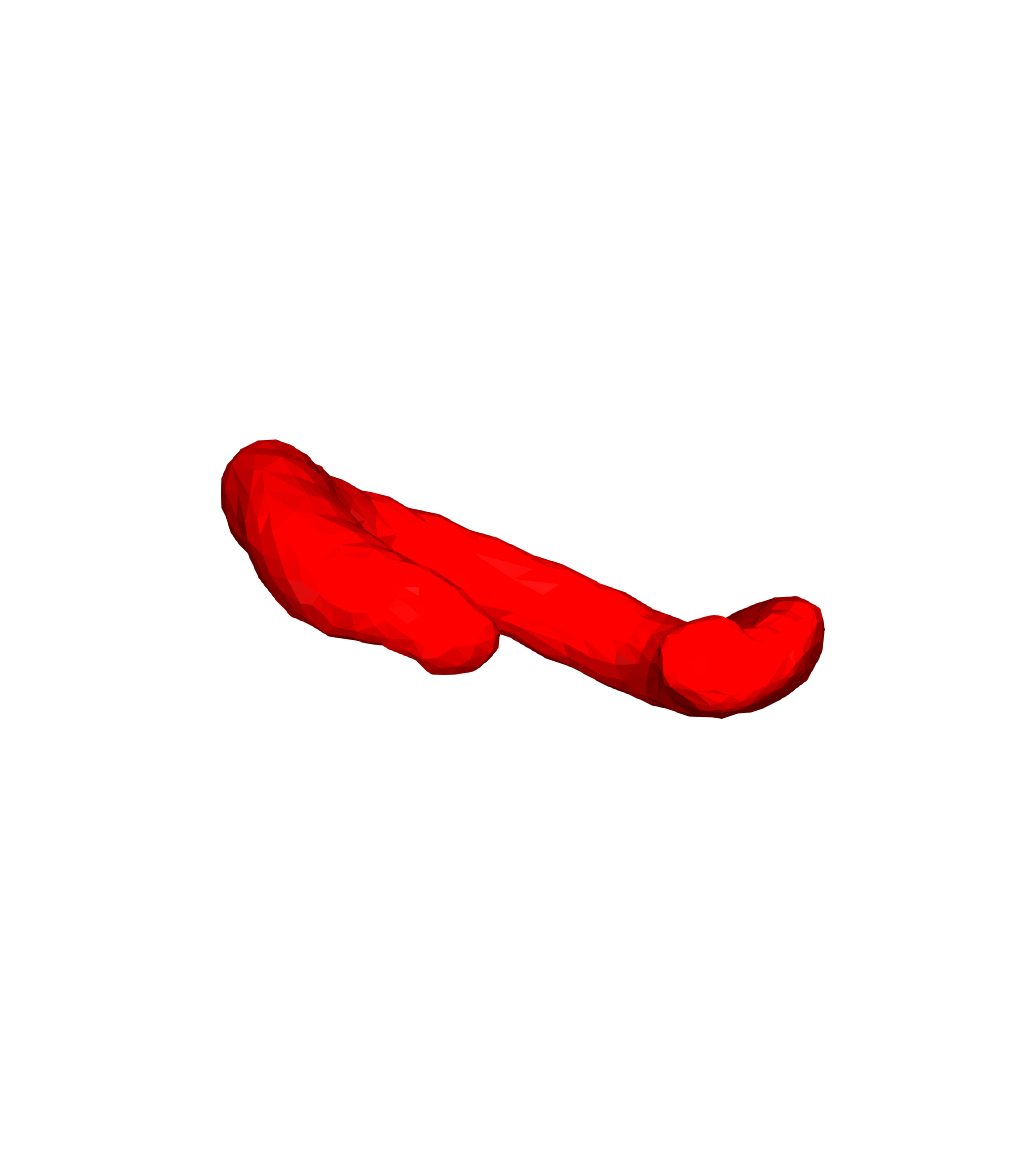}
\includegraphics[trim=5cm 5cm 4.5cm 5cm, clip,width=0.19\textwidth]{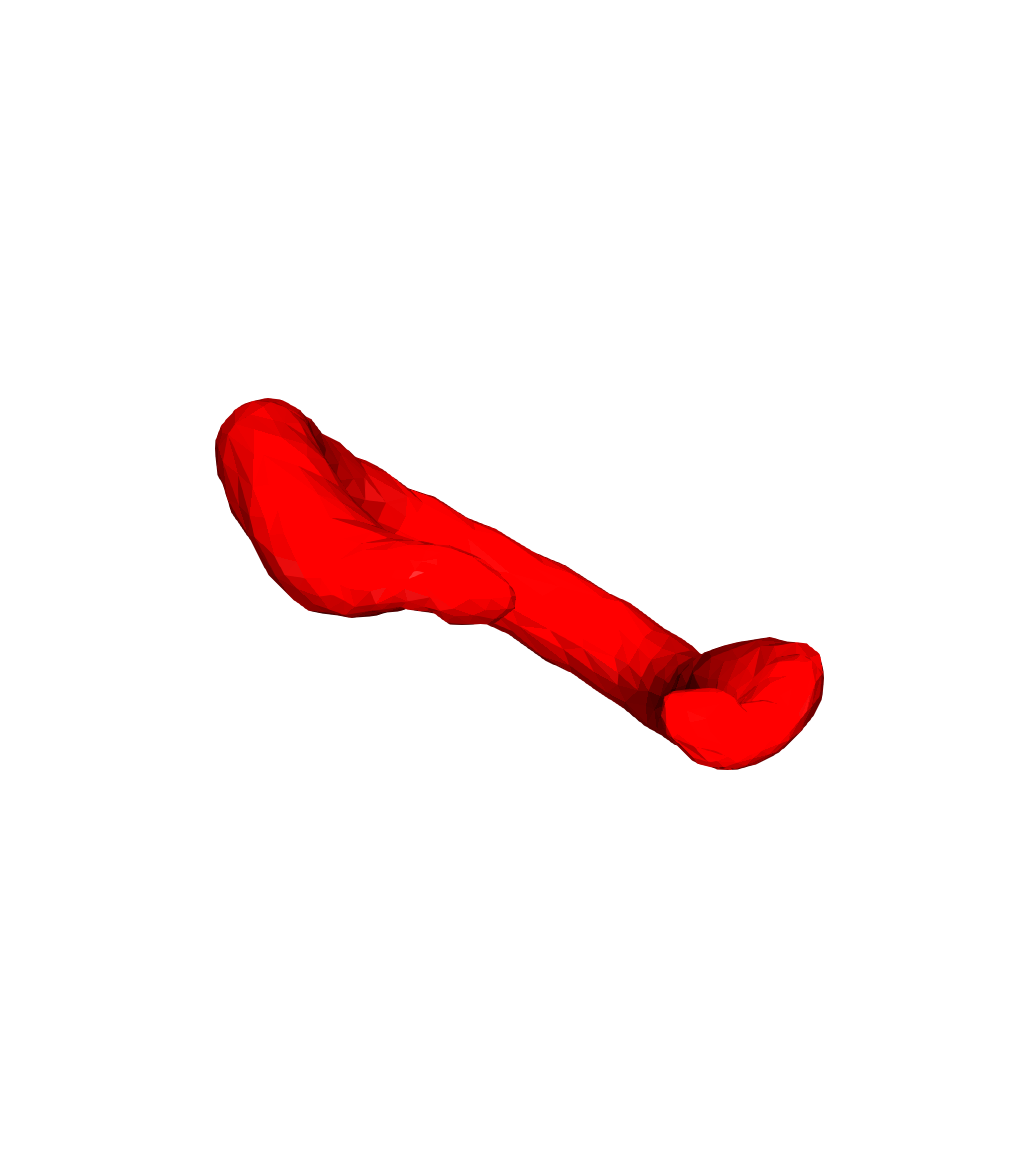}
\includegraphics[trim=5cm 5cm 4.5cm 5cm, clip,width=0.19\textwidth]{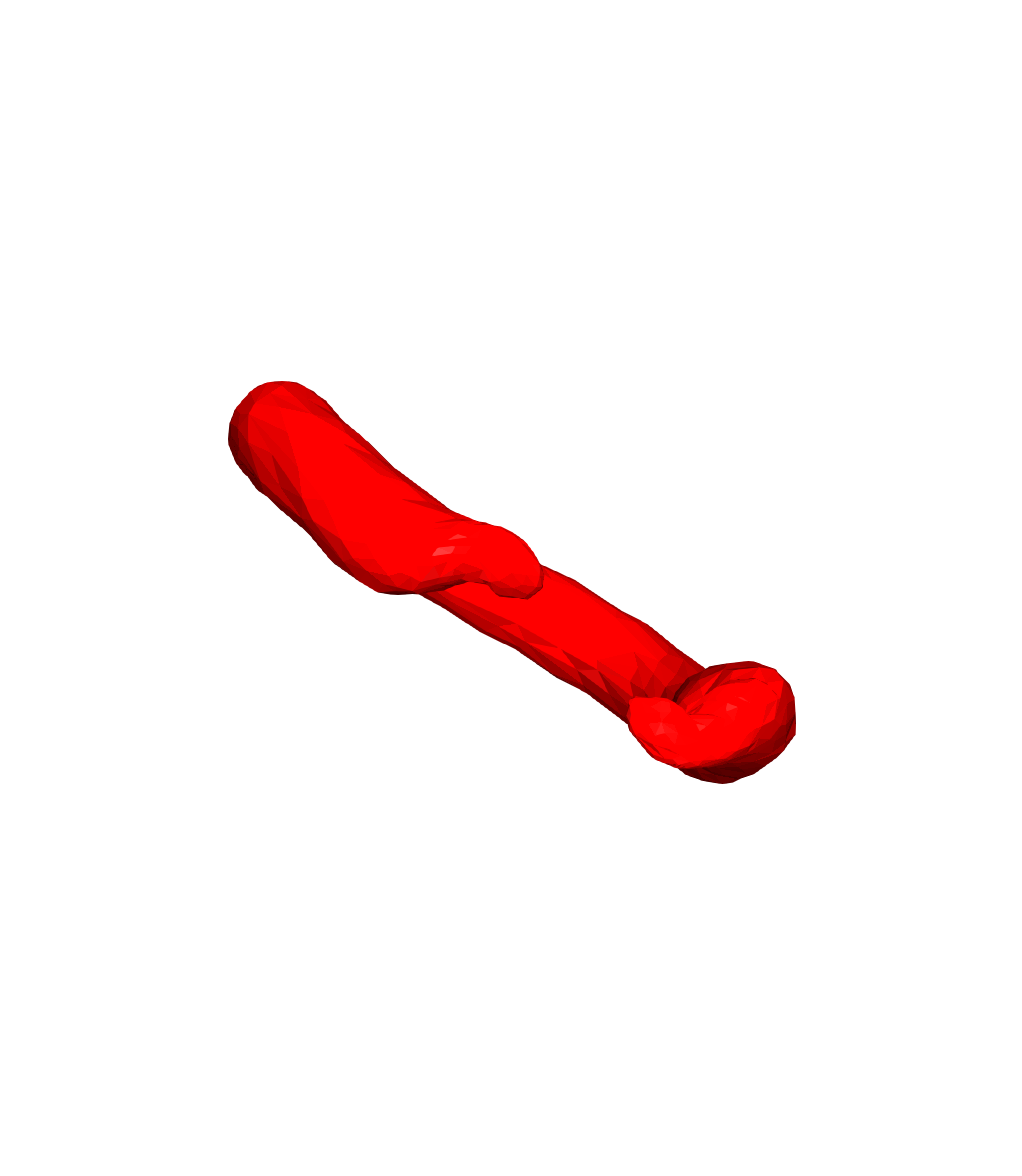}
\includegraphics[trim=5cm 5cm 4.5cm 5cm, clip,width=0.19\textwidth]{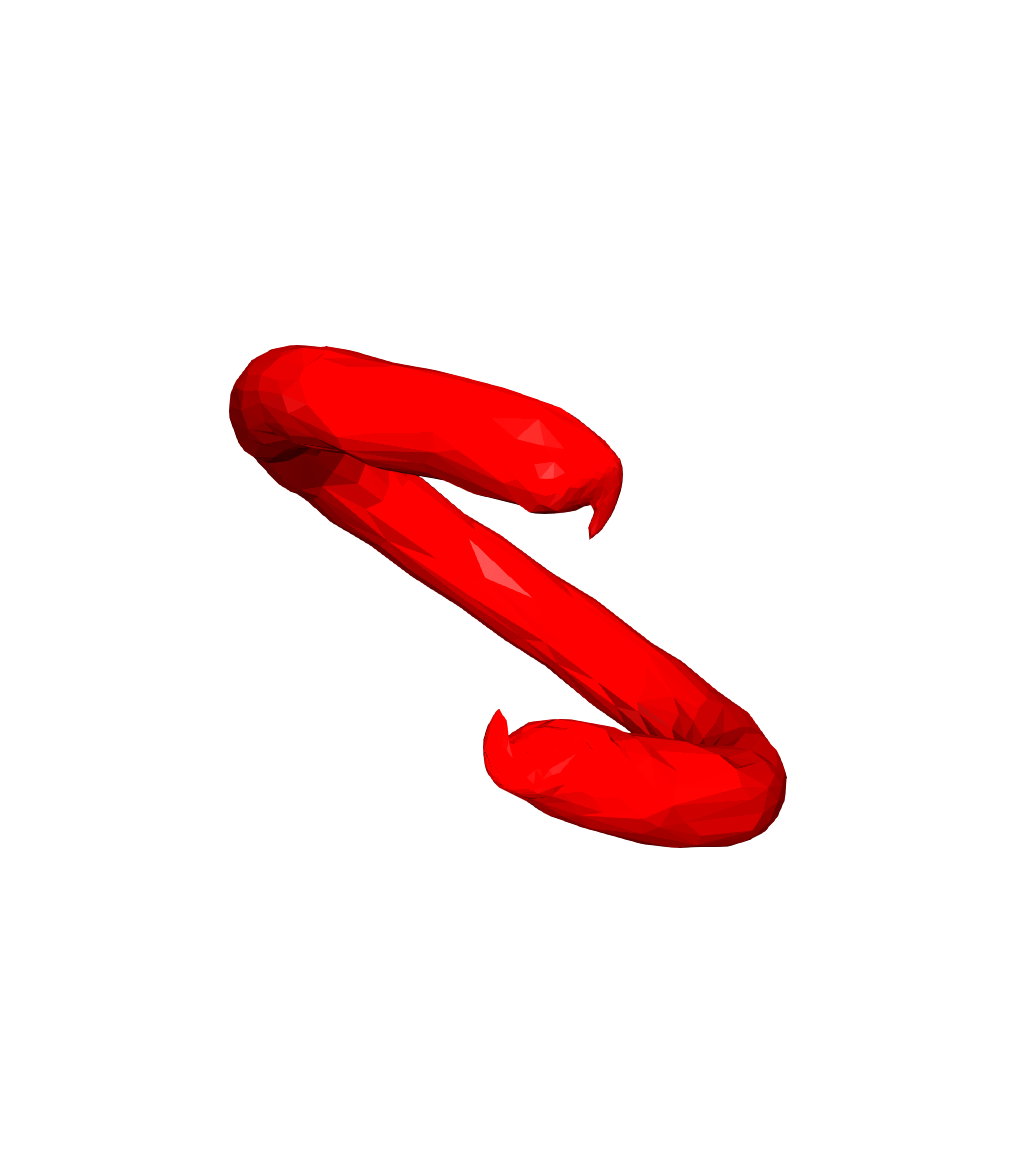}\\

\includegraphics[trim=5cm 5cm 4.5cm 5cm, clip,width=0.19\textwidth]{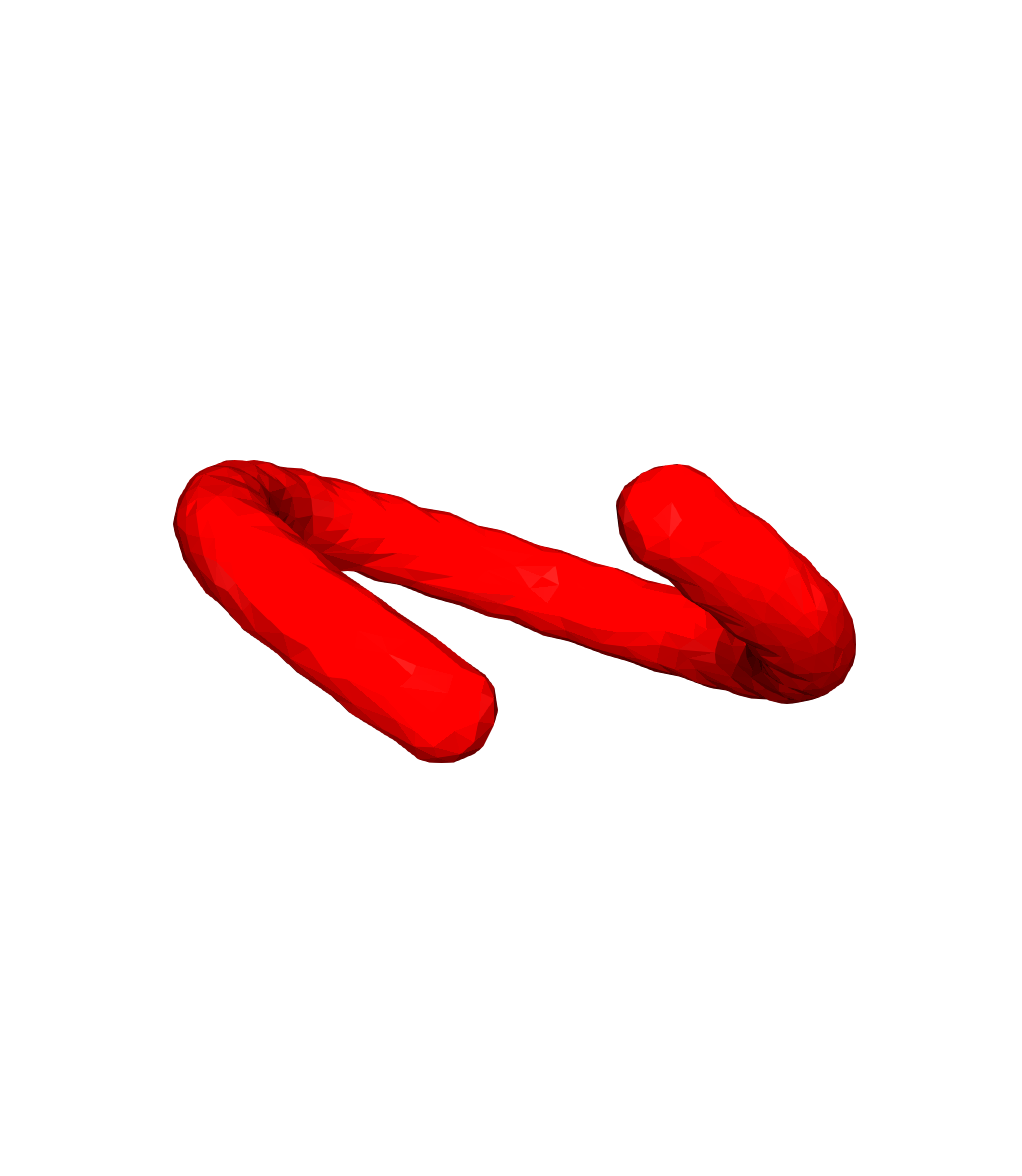}
\includegraphics[trim=5cm 5cm 4.5cm 5cm, clip,width=0.19\textwidth]{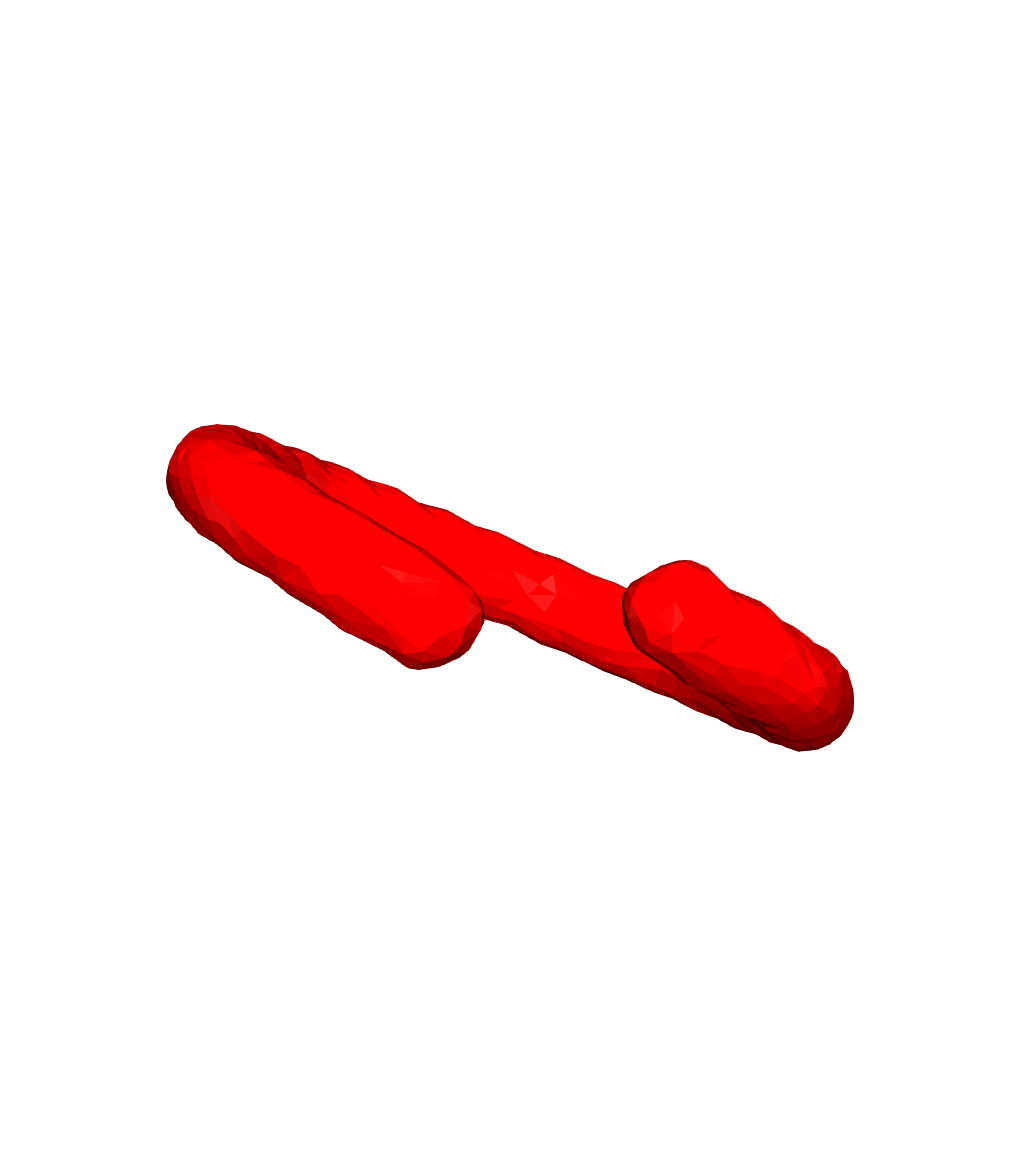}
\includegraphics[trim=5cm 5cm 4.5cm 5cm, clip,width=0.19\textwidth]{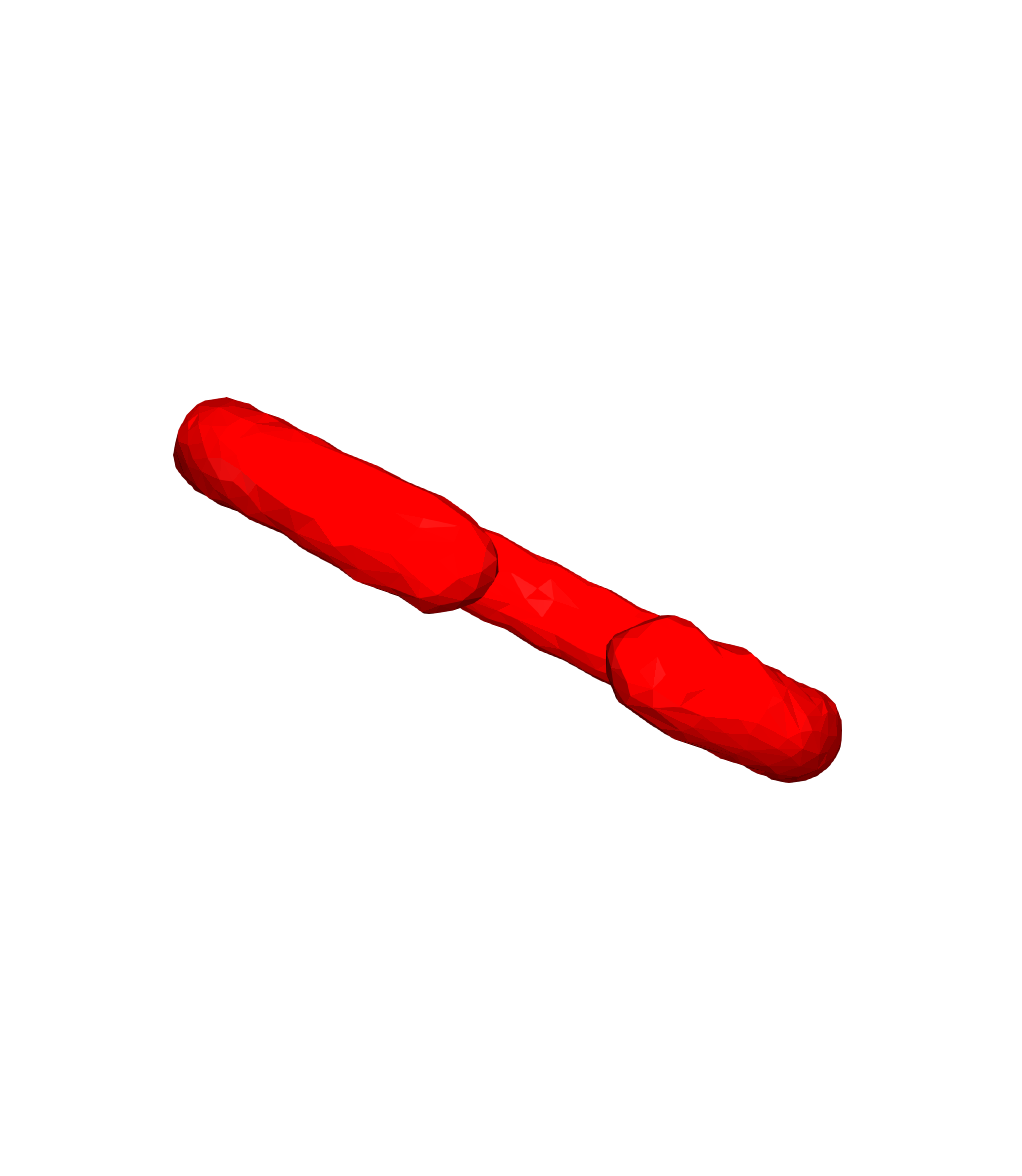}
\includegraphics[trim=5cm 5cm 4.5cm 5cm, clip,width=0.19\textwidth]{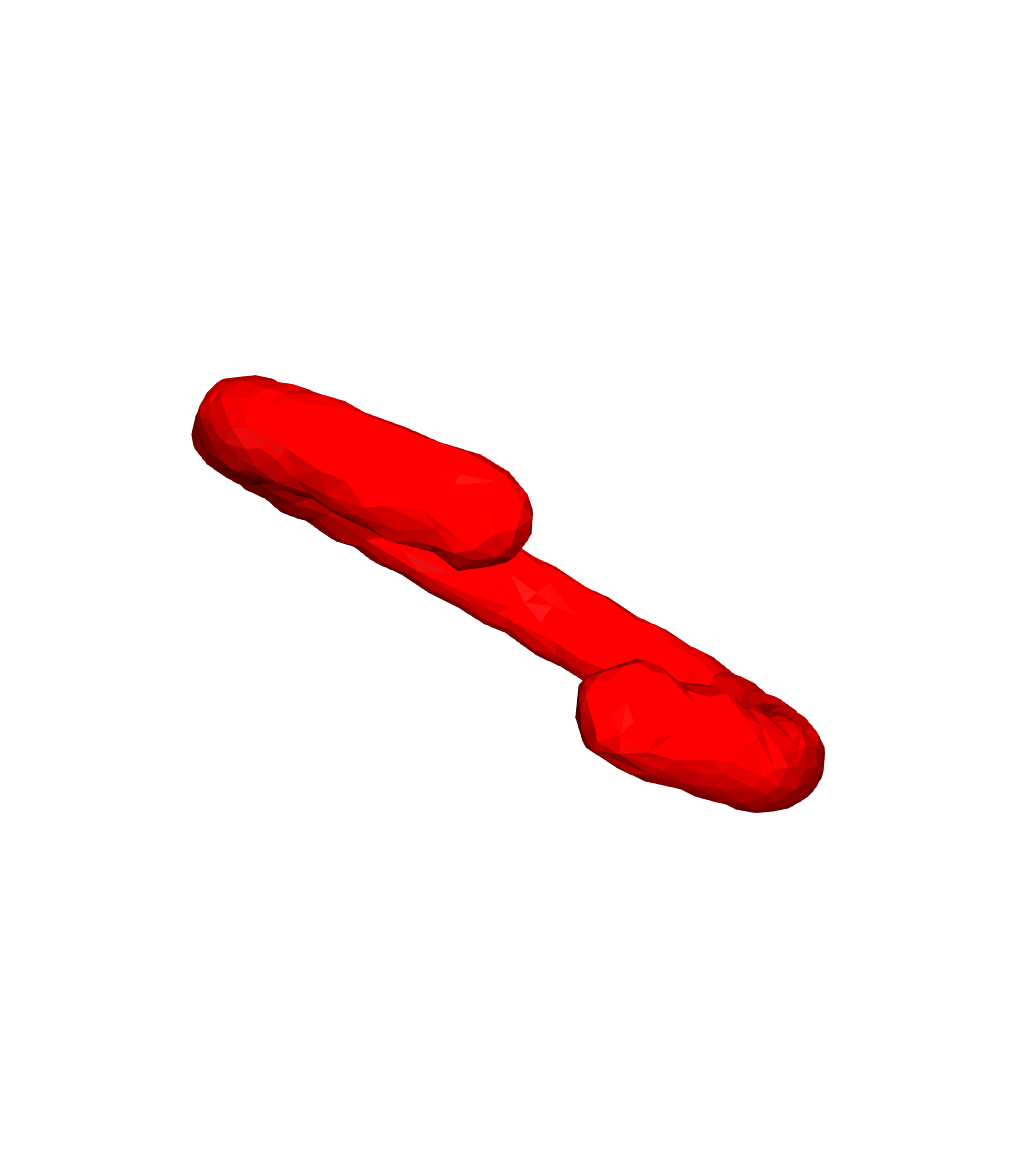}
\includegraphics[trim=5cm 5cm 4.5cm 5cm, clip,width=0.19\textwidth]{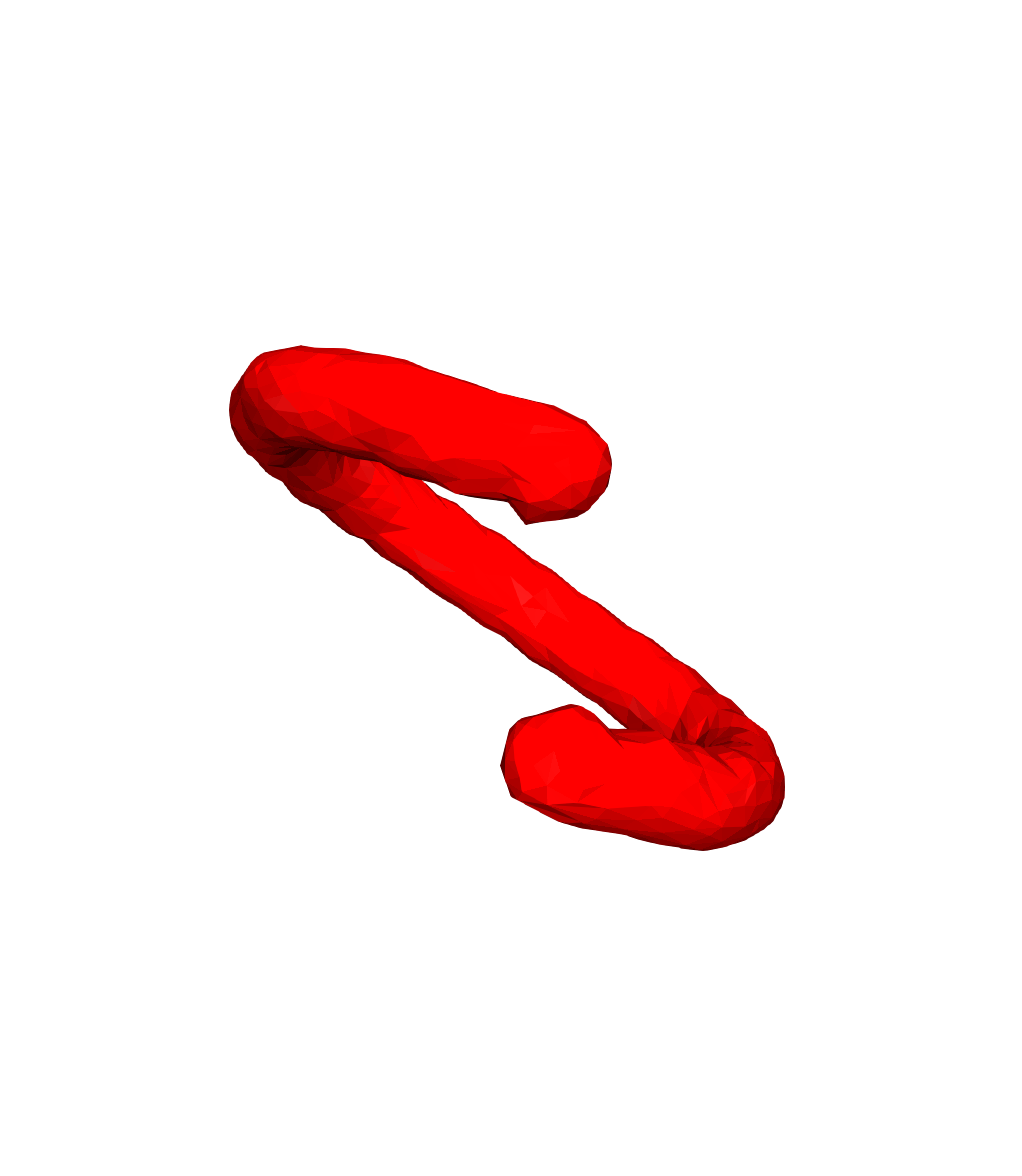}\\

\caption{
\label{fig:outer_hybrid_surfaces} Geodesics between surfaces for outer and hybrid metrics at selected time steps ($t=0, 0.3, 0.5, 0.7, 1$). Top row: outer metric, bottom row: hybrid metric. The hybrid-metric geodesic exhibit significantly less distortion in order to allow the arms of the open ring to pass nearby during the transition from  template to target. }
\end{figure}

The eventual optimization procedure to solve the discrete problem can be then derived from one of the numerical methods for optimal control that were discussed in Section \ref{sec:control}. Figure \ref{fig:outer_hybrid} (second row) shows the geodesic obtained with a hybrid LDDMM + $H^1$ metric on the same curve example as in Figure \ref{fig:geodesic:intrinsic}. In this example, we used an $H^1$ norm penalizing only the tangential part of the derivative of the vector field along the curve (a special case of the two-parameter model introduced in \cite{Mio2004}). Figure \ref{fig:outer_hybrid_surfaces} provides a similar experiment for surfaces, where the hybrid models uses the standard $H^1$ norm discretized according to \eqref{eq:surf.h1}. The surfaces that are matched are deformed incomplete rings, shown in Figure \ref{fig:surf}.

\section{Conclusion}
\label{sec:conclusion}

We have summarized in this chapter some of the most recent developments in metric shape analysis and registration. While this presentation covered a large range of works by our group and others, it certainly does not account for the whole spectrum of literature that involves shape spaces. For example, we only briefly mentioned metamorphosis, which leads to interesting metrics that are still investigated, including recent work such as \cite{hty09,richardson2013computing,Richardson2016,charon2016metamorphoses,effland2018image}. We also  made a short account of the large body of work on the $H^1$ metric on curves or on the LDDMM algorithms and  applications of these two methods, the relevant literature being too large to be listed in this chapter (but the reader may refer to books such as \cite{grenander2007pattern,younes2010shapes,Srivastava2016} and the references they contain). Recent developments on the numerical analysis of this class of problems, include variational time discretization of geodesic calculus \cite{Rumpf2014},  conformal mapping approaches for curve comparison (with the Weil-Petersson metric introduced in \cite{Mumford2006,feiszli2007shape}) and for surface registration \cite{gu2008computational}. One can also cite computational work around the Gromov-Hausdorff (or Gromov-Wasserstein) distance between surfaces \cite{memoli2007use,memoli2011gromov,bronstein2010gromov} or recent methods relaxing of the registration problem through the introduction of functional maps \cite{ovsjanikov2012functional}. \\

The framework discussed in this chapter offers multiple open theoretical or computational problems. In addition to those that were already mentioned in the main text, there is the important issue of generalizing the models beyond the framework of manifolds, including more general objects such as stratified spaces, or trees and graphs, on which very little has been achieved so far \cite{duncan2016elastic}. Also, an important recent and current body of work is focusing on stochastic processes on shape spaces \cite{trouve2012shape,holm2015variational,arnaudon2017stochastic,arnaudon2017geometric,staneva2017learning,sommer2017bridge}, leading to new avenues in the statistical analysis of shapes. The extensive activity within the field is certainly justified by the originality of the mathematical problems that are raised, combined with the diversity and importance of its applications.

\bibliographystyle{abbrv}
\bibliography{references}

\end{document}